\documentclass[10pt,twoside]{article}

\usepackage[numbers]{natbib}
\usepackage[noend, ruled, vlined, linesnumbered]{algorithm2e}
\usepackage[english]{babel}
\usepackage{setspace}

\usepackage{float}
\usepackage{parskip}
\usepackage{booktabs}
\usepackage{caption, subcaption}
\usepackage{array}
\usepackage{amsmath, amssymb, amsthm,bm}
\usepackage[affil-it]{authblk}
\usepackage{booktabs}
\usepackage{hyperref}
\usepackage{fancyhdr}
\usepackage[toc,page]{appendix}

\usepackage{tikz}
\usetikzlibrary{patterns}
\usetikzlibrary{matrix, positioning, fit}
\usetikzlibrary{shapes.multipart}
\usetikzlibrary{arrows.meta}

\newtheorem{definition}{Definition}
\newtheorem{lemma}{Lemma}
\newenvironment{claim}[1]{\par\noindent\underline{Claim:}\space#1}{}

\newcommand{\PreserveBackslash}[1]{\let\temp=\\#1\let\\=\temp}
\newcolumntype{C}[1]{>{\PreserveBackslash\centering}p{#1}}

\newtheorem{ruleEnv}{Theorem}
\newtheorem{Rule}[ruleEnv]{Rule}

\newcommand{\vpi}{\overline{\pi}}
\newcommand{\vrho}{\overline{\rho}}
\newcommand{\vpii}{\vpi^\text{int}}
\newcommand{\vrhoi}{\vrho^\text{int}}
\newcommand{\vrhod}{\vrho^\text{dim}}
\newcommand{\vpid}{\vpi^\text{dim}}
\newcommand{\vlambda}{\overline{\lambda}}
\newcommand{\Ws}{W_{\text{sum}}}
\newcommand{\Rh}{R_\text{rem}}
\newcommand{\Rt}{R^t_\text{max}}
\newcommand{\Rr}{R^r_\text{max}}
\newcommand{\Zl}{\vlambda'}
\newcommand{\Zr}{\vlambda_R}
\newcommand{\Zi}{\vlambda_\mathbb{Z}^*}
\newcommand{\Zint}{\vlambda_\mathbb{Z}}
\newcommand{\Zia}{\vlambda_\mathbb{Z}'}
\newcommand{\Zib}{\vlambda_\mathbb{Z}''}
\newcommand{\Zb}{\overline{\lambda}^\text{inc}}

\definecolor{azul}{RGB}{80,80,160}
\definecolor{mygreen}{RGB}{80,160,80}

\definecolor{mygreen2}{RGB}{213,232,212}
\definecolor{myorange}{RGB}{255,230,204}
\definecolor{myblue2}{RGB}{218,232,252}
\definecolor{myred}{RGB}{248,206,204}
\definecolor{mygrey}{RGB}{235,235,235}

\definecolor{verde2}{RGB}{213,232,212}
\definecolor{laranja}{RGB}{255,230,204}
\definecolor{myblue}{RGB}{90,90,252}
\definecolor{rosa}{RGB}{255, 166, 249}
\definecolor{roxo}{RGB}{142, 128, 255}
\definecolor{vermelho_claro}{RGB}{252, 210, 220}
\definecolor{amarelo}{RGB}{255, 239, 133}
\definecolor{cinza}{RGB}{204, 204, 204}
\definecolor{rosa_claro}{RGB}{255, 209, 253}
\definecolor{azul2}{RGB}{173,216,230}
\definecolor{aquamarine}{RGB}{127,255,212}
\definecolor{verde}{RGB}{152,251,152}
\definecolor{violeta}{RGB}{238,130,238}
\definecolor{verde_escuro}{RGB}{213,255,212}

\makeatletter
\AtBeginDocument{%
  \def\doi#1{\url{https://doi.org/#1}}}
\makeatother

\newcommand\mline[2]{\begin{tabular}[c]{@{}c@{}}#1\\ #2\end{tabular}}

\begin{document}
\onehalfspacing

\def\mytitle{Solving Cutting Stock Problems via an Extended Ryan-Foster Branching Scheme and Fast Column Generation}
\def\myruntitle{Solving CS Problems via an Extended Ryan-Foster BS and Fast Column Generation}

\fancypagestyle{firststyle}
{
   \fancyhf{}
\topskip 30pt\headsep 0pt\headheight 0pt
\renewcommand{\headrulewidth}{0pt}
}

\fancypagestyle{mypagestyle}{%
  \fancyhf{}
  \fancyhead[OR]{R.F.F. da Silva and R.C.S. Schouery}
  \fancyhead[OL]{\thepage}
  \fancyhead[EL]{\myruntitle}
  \fancyhead[ER]{\thepage}
  \renewcommand{\headrulewidth}{.4pt}
}
\pagestyle{mypagestyle}
\date{}

\title{\Large \bf \mytitle}
\author[1]{Renan F. F. da Silva \thanks{Corresponding Author: renan.silva@students.ic.unicamp.br}}
\author[1]{Rafael C. S. Schouery\thanks{ rafael@ic.unicamp.br}}
\affil[1]{Institute of Computing, University of Campinas (Brazil)}
\date{}
\maketitle
\vspace*{-5ex}
\noindent

\abstract{%
We present a branch-cut-and-price framework to solve Cutting Stock Problems with strong relaxations using Set Covering (Packing) Formulations, which are solved by column generation. The main contributions of this paper include an extended Ryan-Foster scheme, which allows us to use this powerful branching scheme even in non-binary problems by using a conflict propagation lemma; a fast column generation process based on a diversification strategy; custom primal heuristics, enabling us to find optimal solutions for several open instances; and a technique to use a smaller feasibility tolerance in floating-point linear programming solvers, combined with numerically safe methods to produce stronger and safer lower bounds.  Additional performance-improving strategies include a technique that controls the height of the branch-and-bound tree; a variable selection algorithm based on branching history; a new set of dual inequalities; insights to obtain a lean model; and the subset-row inequalities. By employing this comprehensive framework, we overcame the current state-of-the-art concerning the following problems: Cutting Stock, Skiving Stock, Ordered Open-End Bin Packing, Class-Constrained Bin Packing, and Identical Parallel Machines Scheduling with Minimum Makespan. Additionally, a new challenging benchmark for Cutting Stock is introduced.
}%

    
\thispagestyle{firststyle}

\allowdisplaybreaks

\section{Introduction}

The main focus of our paper is the Cutting Stock Problem~(CSP), where we have an unlimited number of stock rolls with length $W \in \mathbb{Z}_+$ and a set $I = \{1, \ldots, n\}$ of items, where each item $i \in I$ has a size $w_i \in \mathbb{Z}_+$ and a demand $d_i \in \mathbb{Z}_+$. The objective is to cut the minimum number of stock rolls to satisfy the demands. A special case of the CSP is the Bin Packing Problem~(BPP), where $d_i = 1$, for all $i \in I$.

Successful Integer Linear Programming~(ILP) formulations for these problems include the Set Covering Formulation (SCF) by~\citet{Gilmore_1961} and the arc-flow formulation (AFF) by~\citet{Carvalho_1999}, both with the same relaxation strength. The SCF, also known as a path-flow formulation, has an exponential size, while the AFF has a pseudo-polynomial size.

Let $z_\text{ILP}$ be the optimal solution value and $z_\text{LP}$ be the optimal linear relaxation value of the SCF\@. Instances of the CSP fall into two classes: those with the Integer Round-Up Property (IRUP), which satisfies $z_\text{ILP} = \lceil z_\text{LP} \rceil$, and those without, called non-IRUP instances. The Modified Integer Round-Up Property (MIRUP) conjecture by Scheithauer (1995) suggests that all CSP instances satisfy $z_\text{ILP} \leq \lceil z_\text{LP} \rceil + 1$. Thus, an effective CSP algorithm should handle both instance classes, which requires distinct techniques.

In IRUP instances, finding an optimal solution is typically challenging and demands effective primal heuristics. A straightforward approach involves employing Commercial Integer Programming Solvers, such as Gurobi and CPLEX, known for their strong general primal heuristics. Given the prohibitive number of variables in large instances, even in the AFF, these solvers are often invoked using a restricted set of variables. The work of~\citet{Loti_2022}, based on AFFs, relies exclusively on this as its primary heuristic strategy has achieved great success, which is the state-of-the-art for CSP and other related problems\@.

Yet, this strategy is most effective in AFF-based approaches, where a moderate set of variables allows for the creation of numerous patterns. In contrast, previous works based on SCF employ two custom heuristics capable of generating new columns. The first is a rounding heuristic utilizing the \emph{Sequential Value Correction} method, as employed by~\citet{Belov_2006}. The second is a diving heuristic called \emph{Limited Discrepancy Search}, utilized by~\citet{Wei_2020} for CSP and by~\citet{Pessoa_2020} for the vehicle routing problem (VRP) and its variants, including the BPP (a special case of the Capacitated VRP\@).

For non-IRUP instances, finding an optimal solution is typically straightforward, but establishing a strong lower bound for proving optimality remains challenging. Two primary strategies are used in the literature to address this issue: cutting planes and branching schemes. Regarding cutting planes for the CSP,~\citet{Belov_2006} employs the Chvátal-Gomory Cuts, while~\citet{Wei_2020} and~\citet{Pessoa_2020} use the Subset-Row Inequalities (a subclass of Rank-1 Cuts). As for branching schemes within the SCF, options include branching in variables, utilized by~\citet{Vance_1998} and~\citet{Belov_2006} for the CSP, the Ryan-Foster scheme, adopted by~\citet{Vance_1994, Wei_2020}, and \citet{Baldacci_2024} for the BPP, the implicit flow network scheme, applied by~\citet{Mrad_2013} for the Two-Staged Guillotine Cutting Stock Problem, and a generic scheme used by~\citet{Vanderbeck_2011} for the CSP\@.

Furthermore, several strategies exist in the literature aimed at reducing model size. \citet{Delorme_2020} proposed the Reflect Formulation, an arc-flow approach grounded in the Meet-In-The-Middle principle, significantly reducing arc numbers in some benchmarks by indexing variables using only half of the roll length. \citet{Delorme_2020} and \citet{Loti_2022} utilized variable elimination by reduced cost. Lastly, \citet{Loti_2022} applied waste optimization, employing a volume bound to eliminate arcs from the arc-flow network that cannot contribute to a solution better than the incumbent solution.

Given the substantial size of large CSP instances, the foremost algorithms in the literature, including~\citet{Loti_2022}, predominantly rely on linear relaxation based on Column Generation~(CG). However, in certain instance classes, such as AI and ANI~\citep{Delorme_2016}, CG faces challenges with slow convergence, necessitating specific mitigation strategies. Dual Inequalities, proven effective by~\citet{Amor_2006} and~\citet{Irnich_2016}, offer one approach. Additionally,~\citet{Loti_2022} proposed a multiple pattern generation, which achieves convergence with fewer calls to the pricing algorithm.

This paper introduces a branch-cut-and-price framework tailored to solve problems with strong relaxations using the SCF or the Set Packing Formulation~(SPF). While our primary focus lies on addressing the CSP, we extend our framework to tackle the Skiving Stock Problem~(SSP), Identical Parallel Machines Scheduling with Minimum Makespan~(IPMS), Ordered Open-End Bin Packing Problem~(OOEBPP), and Class-Constrained Bin Packing Problem~(CCBPP). The main contributions of this paper are: an extended Ryan-Foster scheme, which allows us to use this powerful branching scheme even in non-binary problems by using a conflict propagation lemma; a fast column generation process based on a diversification strategy; custom primal heuristics, enabling us to find optimal solutions for several open instances; and a technique to use a smaller feasibility tolerance in Floating-point Linear Programming Solvers~(FLPSs), combined with numerically safe methods to produce stronger and safer lower bounds.  Additional performance-improving strategies include the \emph{splay operation}, a technique that controls the height of the branch-and-bound~(B\&B) tree; a variable selection algorithm based on branching history; a new set of dual inequalities; insights to obtain a lean model; and the subset-row inequalities. By adopting this comprehensive framework, we overcame the current state-of-the-art in the five studied problems, demonstrating the framework's extensibility and its potential applicability to a broader range of problems.

A preliminary version of this paper was published as part of an undergraduate research award competition of a Brazilian conference~\citep{Silva_2024}.

This preliminary version lacks the depth and detail of the current paper, as it is a brief report on the research conducted. Thus, the preliminary version has all proofs omitted, including the proof of the conflict propagation lemma, which is essential for our algorithm's correctness. It also omits important details, such as how to implement the variable selection rule based on branching history, for example.

This paper also presents new contributions that are not included in the preliminary version. It introduces a new set of dual inequalities and a numerically safe method to obtain lower bounds. It also extends the results of the preliminary version to three additional problems (IPMS, OOEBPP, and CCBPP), whereas the preliminary version only considers the CSP and SSP\@. Furthermore, it discusses the necessary adaptations to address these problems and presents detailed computational results for all five problems. Finally, it introduces a new benchmark for the CSP, which may be of great interest for future research and is not included in the preliminary version.

The remaining sections of this paper are structured as follows. Section~\ref{sec::prem} outlines our branch-cut-and-price algorithm, detailing the components including our pricer, cutting planes, branching scheme, and an overview of other key features of our framework. Sections~\ref{sec::conflicts} and~\ref{sec::num_safe} introduce our conflict propagation mechanism and the method for using a smaller tolerance in FLPSs while ensuring numerical safety, respectively. Section~\ref{sec::CRG} develops the column-and-row generation algorithms.
The implementation details of the Relax-and-Fix primal heuristic, along with a variant proposed by us, are explained in Section~\ref{sec::PH}. In Section~\ref{sec::Optimization}, we describe our branching rule, the splay operation, our dual inequalities, the waste optimization, and a variable elimination by reduced cost method, where the latter two are aimed at achieving a leaner model. Furthermore, Section~\ref{sec::CE} presents the computational experiments conducted for the CSP, while Section~\ref{sec::Problems} demonstrates how our framework can be adapted for the other four problems studied. Finally, Section~\ref{sec::Conclusions} draws our conclusions and outlines potential avenues for future research. 

\section{Branch-Cut-and-Price Algorithm\label{sec::prem}}

Our branching scheme is a generalized Ryan-Foster scheme. At each B\&B's node, we maintain a set of items $\overline{I} = \{1, \ldots, \overline{n}\}$, where each item $i \in \overline{I}$ has a demand $d_i > 0$ and a distinct size $w_i > 0$ (i.e., all items of the same size are grouped as copies of the same item), along with a conflict graph $G = (\overline{I}, E)$, in which each edge $\{i, j\} \in E$ indicates that item $i$ cannot be packed with item $j$ and, initially, $E = \emptyset$. The scheme selects a pair of items $(i, j)$. On the left branch, it merges a copy of $i$ with a copy of $j$ into a copy of an item $k$ with size $w_k = w_i + w_j$. In fact, it reduces $d_i$ and $d_j$ by one unit and sets $d_k = 1$, if no item of size $w_k$ exists, or increases $d_k$ by one unit, otherwise. On the right branch, it adds a conflict between all items with sizes $w_i$ and $w_j$. The formal details of this branching scheme are presented in Section~\ref{sec::scheme} and these conflicts are handled by our pricer in Section~\ref{sec::CRG}.

In the context of the CSP, a pattern $P$ represents a feasible way of cutting a stock roll to obtain the demanded items, satisfying the condition $\sum_{i = 1}^{\overline{n}} a_i^{P} w_i \leq W$, where $a_i^P$ denotes the number of copies of item $i$ in~$P$. A pattern is considered valid if for every pair of items $\{i, j\} \in P$ they are not in a conflict, i.e., $\{i, j\} \notin E$. Thereafter, we will simply use the term \emph{pattern} to refer to a valid pattern. A possible formulation for the CSP is the SCF proposed~\citet{Gilmore_1961}, which we describe next:%
\begin{alignat}{3}
  \quad& minimize  \quad  && \displaystyle{\sum_{P \in \mathcal{P}} \lambda_P}\\
              & subject~to \quad       &&\displaystyle{\sum_{P \in \mathcal{P}} a_{i}^P \lambda_P \geq d_i,}    \quad&& \forall i \in \overline{I},\\
              &           \quad       &&\lambda_P \in \mathbb{Z}_+, \quad    && \forall P \in \mathcal{P},
\end{alignat}%
where $\mathcal{P}$ is the set of all patterns, and, for each $P \in \mathcal{P}$, there is an integer variable $\lambda^P$, which represents how many times pattern $P$ is used. The objective is to minimize the number of patterns used.

The formulation (1)--(3) can be enhanced by cutting planes. Let $\mathcal{T} = \{T \subseteq \overline{I}: |T| = 3 \wedge (d_i~=~1, \forall i \in T)\}$ be the set of all triples of items with unitary demands, and $\mathcal{P}(T) = \{P \in \mathcal{T} : |P \cap T| \geq 2\}$ be the subset of patterns with at least two items in $T$. Thus, the weak Subset-Row Inequalities (SRI) of size 3 are defined as
\begin{equation}
    \sum_{P \in \mathcal{P}(T)} \lambda_P \leq 1, T \in \mathcal{T}.\label{eq::cuts}
\end{equation}

These inequalities are well-known for defining facets of the Set Partition Formulation ~\citep{Bulhoes_2018}, which is also a valid formulation for the CSP\@. The presented inequalities are termed \emph{weak} because they are limited to items with unitary demands. SRI and its generalization, known as Rank 1 Chavátal-Gomory inequalities, have seen widespread use in recent decades due to their efficacy in practical problems, especially in BPP~(see, e.g.,~\citet{Wei_2020, Baldacci_2024}) and VRP (see, e.g.,~\citet{Jepsen_2008, Pessoa_2020}).

We refer to (1)--(4) as our Master Problem~(M), with its linear relaxation denoted as LM\@. Due to the potentially large sizes of $\mathcal{P}$ and $\mathcal{T}$, we solve LM using Column-and-Row Generation (CRG). The Restricted Master Problem (RM) is denoted by the model (1)--(4) utilizing only the subsets $\overline{\mathcal{P}} \subseteq \mathcal{P}$ of patterns and $\overline{\mathcal{T}} \subseteq \mathcal{T}$ of cuts, with its linear relaxation referred as RLM\@. The dual of LM is described as follows:
\begin{alignat}{3}
    \quad& maximize  \quad  && \displaystyle{\sum_{i \in \overline{I}} \pi_i} + \displaystyle{\sum_{T \in \mathcal{T}} \rho_T}\\
                & subject~to \quad       &&\displaystyle{\sum_{i \in P} a_{i}^P \pi_i + \sum_{T \in \mathcal{T} : P \in \mathcal{P}(T)} \rho_T \leq 1,}    \quad&& \forall P \in \mathcal{P},\\
                &           \quad       &&\pi_i \in \mathbb{R}^+, \quad    && \forall i \in \overline{I},\\
                &           \quad       &&\rho_T \in \mathbb{R}^-, \quad    && \forall T \in \mathcal{T},
  \end{alignat}%
where $\pi$ and $\rho$ are the dual variables associated with constraints (3) and (4), respectively.

Let $(\vlambda, (\vpi, \vrho))$ be a pair of primal and dual optimal solutions for RLM\@. Our pricing problem seeks a pattern $P \in \mathcal{P}$ with the smallest reduced cost $\overline{c}_P(\vpi, \vrho, b) = b - \sum_{i = 1}^n a_i^P \vpi_i - \sum_{T \in \mathcal{T} : P \in \mathcal{P}(T)} \vrho_T$, where $b = 1$ in this case, though it may vary throughout the paper. In the CRG algorithm, we solve RLM to obtain $(\vlambda, (\vpi, \vrho))$ and execute the pricer to identify a pattern $P$ minimizing $\overline{c}_P(\vpi, \vrho, 1)$. If $\overline{c}_P(\vpi, \vrho, 1) < 0$, then the constraint of $P$  is violated in (5)--(8), consequently, we add $P$ to RLM and re-optimize it. Otherwise, $\overline{c}_P(\vpi, \vrho, 1) \geq 0$ for all $P \in \mathcal{P}$, thus we verify, by enumeration, if there is a cut $T \in \mathcal{T}$ violated by $\vlambda$. If so, we add $T$ to RLM and re-optimize it, and start the process by searching for new violated patterns again. Otherwise, $(\vlambda, (\vpi, \vrho))$ are a pair of primal and dual optimal solutions for LM\@. To speed up its execution, our CRG generates multiple patterns and cutting planes at each step, and it is detailed in Section~\ref{sec::CRG}.

A B\&B that uses a formulation based on CRG is referred to as branch-cut-and-price~(BC\&P). In the following sections, we present our branching scheme in more details and provide an overview of our BC\&P framework, outlining additional features used to enhanced its performance.

\subsection{Extended Ryan-Foster Scheme\label{sec::scheme}}
The overall goal of a branching scheme is to minimize the number of B\&B nodes that need to be evaluated. As SCF and SPF offer strong relaxations for the studied problems, there is almost no integrality gap to close. Thus, the quality of a branching scheme is typically related to its balance of the search space partition and its impact on the pricing substructure. For instance, while branching in variables can result in highly unbalanced partitions and only slightly degrade the pricing substructure, the Ryan-Foster scheme tends to produce more balanced partitions but significantly degrades the pricing structure. Other schemes, such as the implicit flow network scheme~\citep{Mrad_2013} and the generic scheme proposed by~\citet{Vanderbeck_2011}, do not alter the pricing substructure. However, the former does not generate well-balanced partitions, while the latter increases the number of pricing sub-problems in each branch.

In this work, we adopt the Ryan-Foster scheme, and our pricer accommodates this choice, effectively managing the structural degradation caused by branching and cutting planes. To the best of our knowledge, this scheme has been only used for binary matrix models such as the BPP (see, e.g., \citet{Wei_2020,Baldacci_2024}). Hence, we propose an extension for non-binary problems such as the CSP\@.

Let $\vlambda$ be the optimal primal solution for RLM given by the FLPS\@. We define the affinity $\delta_{ij}$ between items $i$ and $j$ as follows: $\delta_{ij} =~\sum_{P \in \mathcal{P} : \{i, j\} \in p} a_i^P a_j^P \vlambda_P$ if $i \neq j$, and $\delta_{ii} = \sum_{P \in \mathcal{P} : i \in p} \frac{a_i^P (a_i^P - 1)}{2} \vlambda_P$, otherwise. Thus, the affinity is the number of times item $i$ is allocated with item $j$ in the same pattern in solution $\vlambda$. For binary problems, \citet{Vance_1994} demonstrate that if $\vlambda$ is fractional, then there exists a pair of items $(i, j)$ with $\delta_{ij} \notin \mathbb{Z}$. 
Consequently, the Ryan-Foster scheme, for binary problems, selects $(i, j)$ with $\delta_{ij} \notin \mathbb{Z}$, branching left with the constraint $\delta_{ij} = 1$ and right with $\delta_{ij} = 0$. Note that, in such problems, $\delta_{ii} = 0$ for all items $i$. Nonetheless, as we consider non-binary problems, our extended scheme presented later can choose the pair $(i, i)$. Additionally, these constraints are typically included implicitly, as adding them explicitly in the primal model can generate dual variables that are complicated to deal with in the pricer. On the left branch, items $i$ and~$j$ are replaced by item $k$, with $w_k = w_i + w_j$, and patterns containing only one of these items are removed from LM\@. On the right branch, patterns including both $i$ and $j$ are removed, introducing a conflict in the pricing problem indicating that $i$ and $j$ cannot be packed together, i.e., adding $\{i, j\}$ to $E$.

Next, we present the extended Ryan-Foster scheme tailored for the CSP\@. Let $\overline{I} = \{1, \ldots, \overline{n}\}$ denote the set of items in the current iteration, where each item has a distinct size (i.e., if there are $d$ items with size $w$, they are represented by a single item $i$ with size $w_i = w$ and demand $d_i = d$). In this scheme, the left branch assumes that at least one more pair of items with sizes $w_i$ and $w_j$ are grouped together, thereby reducing $d_i$ and $d_j$ by one (if $d_i = 0$ or $d_j = 0$ after this operation, we remove $i$ or $j$ from $\overline{I}$), and increasing $d_k$ by one, where $k$ represents the item with size $w_k = w_i + w_j$ (creating it if no such item exists). Conversely, the right branch assumes that any pair of items with sizes $w_i$ and $w_j$ are not grouped together, introducing a conflict $\{i, j\}$ to $E$. This conflict is also propagated to all potential items with sizes $w_i$ and $w_j$ that could be created through merge operations in the descendants of this node. The correctness of this branching scheme is supported by a lemma concerning conflict propagation, which is detailed in Section~\ref{sec::conflicts}.

Finally, the lemma established by \citet{Vance_1994} does not hold for non-binary problems, i.e.,  there are fractional solutions $\vlambda$ for which $\delta_{ij} \in \mathbb{Z}_+$ for all items $i$ and $j$. Although rare, we observe this, for example, when $\vlambda$ contains a pattern $P$ composed of multiple items $i$, resulting in $a_i^P \vlambda_P \in \mathbb{Z}_+$ and $\vlambda_P \notin \mathbb{Z}_+$. Note that adding a merge or a conflict between $i$ and $j$ is a partition of the search space independently of $\delta_{ij}$ to be fractional. The difference is that now the integrality of $\delta$ is not a sufficient condition for $\vlambda$ to be an integer solution. Thus, when there is no fractional $\delta_{ij}$, we select the pattern $P$ with the highest $\vlambda_P - \lfloor \vlambda_P \rfloor$ value and branch into its two items with the highest sum of sizes.

\subsection{Overview\label{sec::overview}}
In this paper, we adopt a pattern-based formulation, such as the SCF, for all the problems under study. Specifically, the formulation utilized for the OOEBPP and CCBPP is exactly the SCF, while we employ the SPF for the SSP\@. As for the IPMS, we leverage its relationship with the CSP, solving it using our CSP solver along with a binary search. Additionally, at the root node, all pricing sub-problems required by the aforementioned formulations can be solved in pseudo-polynomial time by using dynamic programming.

Next, we present an overview of our framework, a branch-cut-and-price~(BC\&P) algorithm based on SCF\@. Our framework's success is based on primal heuristics, fast-increasing lower-bound techniques, and a fast relaxation solver. We outline the pseudocode of our framework in Algorithm~\ref{alg::framework}.

\begin{algorithm}[t!]
    \small
    Compute volume bound

    Run the initial heuristic to get the initial incumbent and columns

    Solve root node with dual inequalities and binary pricing

    Disable dual inequalities and binary pricing

    Solve the root node again

    current node $\gets$ root node

    \While{solution is not optimal}{
        Runs CRF and RF heuristics if their counters are reached

        current node $\gets$ left child from the current node

        \While{current node can be pruned by the numerically safe dual bound}{

            \If{it is possible to remove nodes using splay operation}{
                Remove nodes using the splay operation
            }
            \ElseIf{there is an unexplored right node}{
                current node $\gets$ the deepest unexplored right node}
            \Else{
                break
            }
        }
    }
  \caption{\label{alg::framework}Our framework.}
\end{algorithm}

Our algorithm starts by calculating the volume bound $\lceil (\sum_{i \in I} d_i w_i)/W \rceil$ (or $\lfloor (\sum_{i \in I} d_i w_i)/W \rfloor$ for SSP, as it is a maximization problem). We then solve the root node using new dual inequalities and binary pricing, generating patterns with at most one copy of each item size. These strategies quickly stabilize the dual variable values. In Line 5, we proceed with the CG process without these features, ensuring the algorithm's correctness. Our CG uses a proposed technique to produce a smaller feasibility tolerance in the FLPS, preventing premature stopping, and converging in a few iterations by using a multiple pattern generation with a diversification strategy. The algorithm's correctness is ensured using a pricer with fixed-point numbers and numerically safe dual bounds as also used by~\citet{Baldacci_2024}.

Next, the BC\&P phase begins. At each node, we strengthen the RLM solution from the SCF using weak SRIs with three rows. Our algorithm incorporates custom heuristics, which are activated based on predefined thresholds. As custom heuristics previously used in literature are insufficient to solve the challenging instance class AI of CSP\@, we present a novel approach that combines Relax-and-Fix (RF)~\citep{Belvaux_2000} and a variant proposed by us named Constrained Relax-and-Fix (CRF). In comparative evaluations, our custom heuristics demonstrate competitiveness with Gurobi's general heuristics, as utilized by~\citet{Loti_2022}, the current state-of-the-art for CSP\@.

Subsequently, we adhere to a depth-first order in the BC\&P, prioritizing the left branch. Our branching scheme, the extended Ryan-Foster scheme that we propose, employs a rule based on past decisions and incorporates a conflict propagation strategy among items to break certain symmetries of the problem. 

If the current node can be pruned by the numerically safe dual bound, we employ a new strategy termed \emph{splay operation}, which drops off potentially unfruitful nodes from the BC\&P tree to reduce the execution time. We reevaluate the current node if this strategy results in the removal of at least one node. Otherwise, we proceed to explore the next unexplored right node. Furthermore, we terminate the algorithm as soon as we detect an optimal solution, which may already occur in Line 2.

Finally, we implement a lean RLM to enhance our algorithm's efficiency by implementing two strategies to reduce its size. The first strategy involves limiting the maximum waste of a pattern while solving both the root node and subsequent BC\&P nodes. This technique was employed by~\citet{Loti_2022}, but it was applied only after the root node was solved (not during its solution).  The second strategy entails a routine that relies on reduced costs to temporarily eliminate patterns that cannot belong to an optimal solution. These several techniques give a fast relaxation solver and a fast lower bound increase, enabling the efficient resolution of non-IRUP CSP instances, such as those in the ANI class.

\section{Conflict Propagation\label{sec::conflicts}}
We define a \emph{superitem} $i$ as an item composed of one or more items of $I$, with $0 \leq w_i \leq W$. As previously defined, our set $\overline{I}$ represents a set of superitems where if $i \neq j$, for $i, j \in \overline{I}$, then $w_i \neq w_j$. It is worth noting that the composition of a superitem is inconsequential in this definition of $\overline{I}$; if there are $k$ superitems with the same size $w$, we consolidate them into one item $i$ with demand $d_i = k$ and size $w_i = w$. Henceforth, a superitem is simply referred to as an item.  Additionally, we denote the set of items that has a conflict with an item $i$ in $G = (\overline{I}, E)$ by $E_i$, where $G$ is the conflict graph of $\overline{I}$. While aggregating items of the same size proves advantageous, initially, it may seem necessary to maintain $d_i$ distinct conflict lists $E_i$ for each item $i \in \overline{I}$, i.e., it seems that $G$ should have $\sum_{i \in \overline{I}} d_i$ vertices, not $|\overline{I}|$, similarly to a BPP's instance\@. 

Definition~\ref{definition::equivalance} outlines the concept of equivalence classes of solutions for the CSP\@. Equivalent solutions are common, especially when the demands are non-unitary. With this in mind, Definition~\ref{definition::scheme} presents a scheme that maintains only one conflict list per item size by propagating conflict among items after branching. This conflict propagation avoids analyzing more than one solution of each equivalence class, and its correctness is asserted by Lemma~\ref{lemma::propagation}, with the detailed proof presented in Appendix~\ref{app::propagation}.

\begin{definition}\label{definition::equivalance}
    Given a solution $\Zia$, note that swapping an item of size $w$ with a set of items $S$, where $\sum_{i \in S} w_i = w$, does not affect feasibility or solution value. Therefore, we define two solutions $\Zia$ and $\Zib$ as equivalent if we can transform $\Zia$ into $\Zib$ using a sequence of these swap operations.
\end{definition}

\begin{definition}\label{definition::scheme}
    Let $(i, j)$ be the pair of items chosen for branching at the current node, and let $k$ be an item such that $w_k = w_i + w_j$. Assuming that all items with the same size share the same conflict list, denote $E_h$ and $E_h'$ as the conflict lists for all items with size $w_h$ before and after branching at the current node, respectively. We define $E_k' = E_k \cup E_i \cup E_j$ on the left branch, and $E_i' = E_i \cup \{j\}$ and $E_j' = E_j \cup \{i\}$ on the right branch.
\end{definition}

\begin{lemma}\label{lemma::propagation}
    If the B\&B tree is traversed using a depth-first search favoring the left branch, then the scheme outlined in Definition~\ref{definition::scheme} is correct, as all discarded solutions are equivalent to those already analyzed by the algorithm.
\end{lemma}

\section{Numeric Safety and Smaller Tolerance\label{sec::num_safe}}

The fastest linear programming solvers such as Gurobi are floating-point based and use a tolerance parameter~$\epsilon$ to avoid convergence problems, which is necessary due to the floating-point arithmetic. In this context, a dual solution $(\vpi, \vrho)$ is considered feasible if it satisfies ${\sum_{i \in P} a_{i}^P \vpi_i + \sum_{T \in \mathcal{T} : P \in \mathcal{P}(T)} \vrho_T \leq 1 + \epsilon}$ for each constraint $P$, i.e., the solver accepts an infeasible solution as feasible if it violates each constraint by a value less than $\epsilon$. Consequently, the value of $(\vpi, \vrho)$ may not be a valid lower bound. Moreover, it is only beneficial for the pricer to find patterns with reduced cost $\overline{c}_P(\vpi, \vrho, 1) <~-\epsilon$, as other patterns may not be included in the RLM basis by the solver. Finally, if the pricer relies on floating-point numbers, a solution may appear feasible due to numerical cancellation, even when it is, in fact, infeasible. Hereafter, we call \emph{floating-optimal} the primal and dual solutions claimed to be optimal by the FLPS\@.

Both~\citet{Baldacci_2024} and~\citet{Loti_2022} address these issues to achieve a safe dual bound. Specifically,~\citet{Baldacci_2024} employ a pricer with fixed-point numbers to produce a numerically safe algorithm as well, contrasting with~\citet{Loti_2022}, whose approach remains susceptible to floating-point issues. Both studies introduce techniques applicable to formulations based on column-and-row generation, expanding upon the method used by \citet{Held_2012}, which is limited to column generation.

Given $(\vpi, \vrho)$ a floating-optimal dual solution, and two integers $K$ and $M$, such that ${K = 2^{k} \geq M = 2^{m}}$ for some $k, m \in \mathbb{Z}_+$, and $\frac{1}{M} \geq \epsilon$. We denote the scaled integer dual solution by $(\vpii, \vrhoi)$, with $\vpii_i = \lfloor K \vpi_i \rfloor$ for all $i \in \overline{I}$ and $\vrhoi_T = \lfloor K \vrho_T \rfloor$, for all $T \in \overline{T}$, and we denote the (fixed-point) diminished dual solution by $(\vpid, \vrhod)$,  with $\vpid_i = \frac{1}{K} \vpii_i $ for all $i \in \overline{I}$ and $\vrhod_T = \frac{1}{K} \vrhoi_T $, for all $T \in \overline{T}$. \citet{Baldacci_2024} proves that $$z^\text{safe} = \frac{z((\vpid, \vrhod))}{1 - \overline{c}(\vpid, \vrhod, 1)}$$ is a numerically safe lower bound if the minimum reduced cost $\overline{c}(\vpid, \vrhod, 1)$ is computed without floating-point errors. This can be achieved by utilizing the objective $\overline{c}_P(\vpii, \vrhoi, K)$ in the pricer, ensuring $(\vpii, \vrhoi)$ can be represented by 64-bit integers. This approach guarantees exact computations (free from floating-point errors) in the pricer if we guarantee that there is no under or overflow. Additionally, $z^\text{safe}$ is a valid lower bound because each constraint in $(\vpid, \vrhod)$ is violated by at most $-\overline{c}(\vpid, \vrhod, 1)$. By dividing all dual variables of $(\vpid, \vrhod)$ by ${1 - \overline{c}(\vpid, \vrhod, 1)}$, we obtain a feasible dual solution with the value $z^\text{safe}$.

Note that, if there is a pattern $P$ with $\overline{c}_P(\vpii, \vrhoi, K) < \frac{-K}{M}$, then $\overline{c}_P(\vpid, \vrhod, 1) < \frac{-1}{M} \leq -\epsilon$. Since all dual values were not increased in the diminished dual solution, it follows that $\overline{c}_P(\vpi, \vrho, 1) \leq \overline{c}_P(\vpid, \vrhod, 1)$, and, thus, $P$ will be pivoted when added to RLM\@. Therefore, our CRG process stops when ${\overline{c}_P(\vpii, \vrhoi, K) \geq \frac{-K}{M}}$ for all patterns $P \in \mathcal{P}$. At this point, we employ an alternative pricer, explained in Section~\ref{sec::recover}, to determine a lower bound~$B$ for the minimum reduced cost $\overline{c}(\vpii, \vrhoi, K)$, and we use $B$ to compute $z^\text{safe}$ (which, by the result of~\citet{Baldacci_2024}, is also a valid and safe bound).

In our experiments, we set $K = 2^{49} \approx 5 \cdot 10^{15}$ and could choose $M = 2^{29}$, given that the smallest tolerance $\epsilon$ supported by the Gurobi (the FLPS used by us) is equal to $10^{-9}$ and $\frac{-1}{M} \approx -1.8 \cdot 10^{-9} < -10^{-9}$. However, we use a larger $M$, as explained later. To prevent under or overflow in the pricer, the summation $\sum_{i \in \overline{I}}d_i \cdot \vpii_i$ must not exceed $\lfloor \frac{2^{63} - 1}{2^{49}} \rfloor = 16383$, and $\sum_{T \in \overline{\mathcal{T}}} \vrhoi_T$ should be greater than or equal to $\frac{-2^{63}}{2^{49}} = -16384$, which remains true in our experiments across the studied instances. If this condition does not hold, we could compute $(\vpii, \vrhoi)$ using a smaller value for $K$. 

A noteworthy point is that, in our framework, the FLPS works as an oracle, providing a good guess for a dual solution $(\vpi, \vrho)$ that is likely to have little or no violations. Regardless of the quality of the FLPS's guess or any issues that the FLPS presents, we can compute the valid lower bound $z^\text{safe}$ at any time. Despite the strategies mentioned above, some numerical issues still impact the FLPS's ability to provide better guesses.

The first issue is that certain CSP benchmarks have patterns with reduced costs in the range $(-10^{-9}, 0)$ and, thus, the minimum tolerance $\epsilon$ supported by Gurobi is not enough to pivot these patterns. For instance, in the AI and ANI instances with 600 items or more, \citet{Pessoa_2020} noted that neglecting to pivot these patterns to the simplex basis results in a weaker lower bound. Based on our experiments, we assert that most instances from the AI and ANI classes pose significant challenges due to the premature termination of CRG caused by this restricted tolerance.

The second issue is related to floating-point problems, where a pattern $P \in \overline{\mathcal{P}}$ with $\overline{c}_P(\vpi, \vrho, 1) < -\epsilon$ fails to pivot into the simplex base due to a reduced cost erroneous computed by the solver. While this scenario is rare, we have observed instances in AI and ANI with 800 items or more where optimization results were compromised by this issue. It requires attention because our pricer could regenerate these columns and cause an infinite loop.

With this in mind, we propose the following strategy based on a tolerance $\epsilon' = 2.5 \cdot 10^{-12}$ and by taking the parameter $M = 2^{38}$, where $\frac{1}{M} \approx 3.6 \cdot 10^{-12}$.
The tolerance $\epsilon'$ addresses the first numerical issue by aiming to enforce Gurobi to operate within this tolerance. While not explicitly supported by the solver, we can implicitly achieve this by multiplying the objective function by a factor $C = 400$ and disabling Gurobi's auto-scaling feature (Scale Flag). This adjustment ensures that Gurobi's solved model has reduced costs scaled by $C$, making the solver's absolute tolerance of $10^{-9}$ equivalent to $2.5 \cdot 10^{-12}$ in the original model. We choose this value of $C$ as it appears sufficient to pivot all columns with negative reduced costs in most cases, halting the CRG when the dual solution is feasible under a despicable tolerance, without exacerbating the second numerical issue significantly.

Concerning the second issue, if our algorithm detects it in $(\vpi, \vrho)$, we discard the current simplex basis and obtain a new one by solving RLM using the barrier method, which is faster than the simplex algorithm for solving from scratch. Using the barrier method, we never detected the second numerical issue on $(\vpi, \vrho)$. However, if it does arise, we could halt the CRG process and calculate a safe bound $z^\text{safe}$ with the current smallest reduced cost. Additionally, we have mitigated occurrences of the second numerical issue by setting Gurobi's parameter \textit{Numeric Focus} to 2, prompting the solver to exercise more caution during numerical manipulation. To detect if $(\vpi, \vrho)$ has the second issue, we assess it under $(\vpii, \vrhoi)$ by verifying if $\overline{c}_P(\vpii, \vrhoi, K) < \frac{-K}{M}$ for some pattern $P \in \overline{\mathcal{P}}$, i.e., for the patterns in RLM. Even though this could lead to false positives and false negatives, this is manageable.
A false positive results in running the barrier method unnecessarily, while a false negative prevents the pricer from generating this column, thereby avoiding potential loops. Furthermore, for a false positive to occur, conversion errors from $(\vpi , \vrho)$ to $(\vpii, \vrhoi)$ must generate an accumulated error greater than $(\frac{1}{M} - \epsilon ') \cdot K \approx 640$ in the reduced cost of a pattern, which is not expected.

Similar to the algorithm proposed by~\citet{Baldacci_2024}, ours is also numerically safe, not being subject to producing an incorrect answer due to floating-point arithmetic. This sets our framework apart from other state-of-the-art algorithms~(\citet{Pessoa_2020,Loti_2022}), which rely on commercial integer solvers in sub-routines that may produce inaccurate results in instances prone to numerical issues (such as AI and ANI with 801 items or more).

\section{Column-and-Row Generation\label{sec::CRG}}
Next, we present our pricer and our SRI enumeration process.

\subsection{Pricing\label{sec::pricing}}

Having obtained a floating-optimal dual solution $(\vpi$, $\vrho)$ for RLM, we derive a scaled integer dual solution $(\vpii, \vrhoi)$ and a diminished dual solution $(\vpid, \vrhod)$ using the parameter $K$ and $M$, as defined in Section~\ref{sec::num_safe}. Let $\overline{I} = \{1, \ldots, \overline{n}\}$ denote our set of items, where each item $i \in \overline{I}$ possesses a demand $d_i$ and a size $w_i$. Moreover, consider $G = (\overline{I}, E)$ as a conflict graph, where each edge $(i, j) \in E$ represents a conflict between items $i$ and $j$ introduced by our branching scheme, and $\overline{\mathcal{T}}$ represents the set of active SRIs.

Following the approaches of~\cite{Wei_2020} and~\citet{Baldacci_2024} for the BPP, we can formulate the optimization problem, known as the Knapsack Problem with Triples and Conflicts (KP-T-C), that aims to find a pattern $P$ with the minimum reduced cost $\overline{c}_P(\vpii, \vrhoi, K)$ as follows.
\begin{alignat}{3}
    \text{(KP-T-C)}\quad c(\vpii, \vrhoi, K) = &   \min \quad && \displaystyle{ K - \sum_{i \in \overline{I}} \vpii_i a_i + \sum_{T \in \overline{\mathcal{T}}}} \vrhoi_T b_T \label{form::begin-kp}\\
    &s. t. &&\displaystyle{\sum_{i \in \overline{I}} w_i a_i \leq W,} \quad&&\\
    &&& a_i + a_j \leq 1, \quad && \forall (i, j) \in E, \\
    &&& a_i + a_j \leq b_T + 1, \quad && \forall T \in \overline{\mathcal{T}},  i, j \in T, i \neq j\\
    &&&0 \leq a_i \leq d_i, \text{ with } a_i \in \mathbb{Z}_+, \quad && \forall i \in \overline{I},\\
    &&&b_T \in \{0, 1\}, && \forall T \in \overline{\mathcal{T}} \label{form::end-kp},
\end{alignat}
where $a_i$ represents the quantity of item $i \in \overline{I}$ in the solution, and $b_T$ is a binary variable ensuring that if at least two items from an SRI $T \in \overline{\mathcal{T}}$ are included in the solution, then the corresponding dual variable of $T$ is accounted for in the solution value. Notably, the optimal solution vector $a^*$ defines a pattern $P$ with the minimum reduced cost $\overline{c}_P(\vpii, \vrhoi, K)$. A pattern $P$ is considered violated if $\overline{c}_P(\vpii, \vrhoi, K) < \frac{-K}{M}$ as discussed in Section~\ref{sec::num_safe}. Furthermore, in that same section, it is stated that if no patterns are violated, then $(\vpid, \vrhod)$ is a numerically safe dual feasible solution to LM\@.

As noted by~\citet{Baldacci_2024}, relying on state-of-the-art MILP solvers to tackle this problem is ineffective, especially in our scenario where we need to solve numerous KP-T-C instances due to our CRG-based formulation. Additionally, these solvers do not guarantee numerical safety when dealing with floating-point numbers, a critical requirement for us. We propose a branch-and-bound algorithm that employs the classical, matrix-based, dynamic programming (DP) for the KP as a lower bound and identifies multiple violated patterns if they exist. Another approach to solving the KP-T-C and its generalizations is the \emph{Label-Setting Algorithm}~(LSA), widely used in literature for problems like the BPP \citep{Wei_2020,Baldacci_2024} and VRPs~\citep{Pessoa_2020}. The LSA utilized by \citet{Wei_2020,Baldacci_2024} also returns multiple patterns, focusing on choosing those with the smallest reduced cost. While our algorithm shares similarities with LSA, we prioritize diversification in pattern selection and do not utilize dominance rules.


\subsection{Motivation for diversification}
As highlighted by some authors~\citep{Amor_2006,Irnich_2016}, the primal space of LM for CSP and BPP can be highly degenerated, resulting in slow convergence in CG-based formulations (which is naturally extended to CRG-based formulations). One potentially useful strategy to guide a more stable CG is to generate multiple patterns in the pricing problem, which is a technique already employed by~\citet{Wei_2020,Loti_2022,Baldacci_2024}.  Nevertheless, this approach involves a trade-off. On the one hand, the number of times the pricer is executed tends to be lower, thus the overall runtime can decrease. Conversely, it also increases the likelihood of adding redundant patterns, potentially inflating the cost of solving intermediate linear programs and leading to longer execution times.

\citet{Wei_2020} and~\citet{Baldacci_2024} return at most ten patterns by iteration, choosing those with the smallest reduced cost among those found.
This restriction probably arises because a higher number would result in an excessive selection of redundant patterns. This claim is supported by the observation that the set of patterns $\overline{\mathcal{P}}''$ with the smallest reduced costs may involve the same subset of items $\overline{I}'$. Consequently, perhaps after adding a small subset of $\overline{\mathcal{P}}''$ to RLM, no patterns in this set remain violated. Conversely, there might be another subset of items $\overline{I}''$ in the current dual solution belonging to violated patterns (although with reduced costs higher   than those in $\overline{\mathcal{P}}''$),  which could still be violated even after adding $\overline{\mathcal{P}}''$.

In summary, we believe that the key to achieving a beneficial trade-off lies in diversification, choosing a set of patterns with as many different items as possible. A diversification strategy has already been successfully implemented by~\citet{Loti_2022}, who used bidirectional DP to generate the most violated pattern (i.e., the one with minimum reduced cost) containing each item. However, since~\citet{Loti_2022} employed a robust branching scheme and did not utilize SRIs, it is not immediately clear how to adapt their strategy efficiently to our algorithm. Consequently, our approach takes a different path, focusing on generating patterns with large items and incorporating a degree of diversification based on the intuition that the slow convergence is primarily due to the difficulty in finding good patterns containing these items.

\subsection{Multiple Patterns Generation Algorithm\label{sec::recover}}
Next, we present our pricer called the Multiple Pattern Generation Algorithm, which focuses on generating patterns with large items and with a degree of diversification. First, note that, in the pricer, we can set $d_i$ as the minimum of its original demand $d_i$ and $\left\lfloor \frac{W}{w_i} \right\rfloor$. Additionally, if an item $i$ conflicts with itself (i.e., $(i, i) \in E$), we set its demand $d_i$ to 1. 

In this section, consider that $\overline{I}$ is a sequence with $d_i$ copies of each item $i$ and indexed from 1 to $|\overline{I}|$. We partition $\overline{I}$ into three subsequences: $\overline{I}_1$, containing items belonging to any $T \in \overline{\mathcal{T}}$ with $\vrho < 0$; $\overline{I}_2$, containing items that conflict with another item; and $\overline{I}_3$, containing the remaining items. Consider that each of these subsequences is arranged in non-increasing order of size and that $\overline{I}$ is equal to the concatenation of $\overline{I}_3, \overline{I}_2$ and $\overline{I}_1$, in this order.

We can obtain a relaxation of (\ref{form::begin-kp})--(\ref{form::end-kp}) by setting $\overline{\mathcal{T}} = E = \emptyset$. This relaxation can be solved in time $\Theta(|\overline{I}| W)$ by using the classical DP for the KP\@. This DP gives us a matrix $dp$, where each cell $dp(i, r)$ has the optimal solution value considering only the first $i$ items of $\overline{I}$ and bin capacity $r$. In other words, $dp(i, r) = K - \sum_{i \in P} \vpii a_i^{P}$, where $P \in \overline{\mathcal{P}}$ is a pattern that minimizes $dp(i, r)$.

We can, then, solve (\ref{form::begin-kp})--(\ref{form::end-kp}) by performing a B\&B employing the matrix $dp$ as a lower bound. In this B\&B, we initialize at state $(|\overline{I}|, W)$ and construct a partial pattern $\overline{P}$, initially empty. We maintain a pool of violated patterns throughout the process. In each state $(i, r)$, the left branch includes item $\overline{I}[i]$, while the right branch does not. The left branch is valid only if $\overline{I}[i]$ is compatible with all items in $\overline{P}$ and $w_{\overline{I}[i]} \leq r$. We prioritize exploring the left branch, and we make a recursive call to a state $f(i, r)$ only if ${f(i, r) - v(\overline{P}) < \frac{-K}{M}}$, where $v(\overline{P}) = \sum_{i \in \overline{P}} \vpii_i + \sum_{T \in \overline{\mathcal{T}}(\overline{P})} \vrhoi_T $. If this condition is not met, the state cannot produce a violated pattern\@. The base cases occur when $i = 0$ or $r = 0$, indicating that $\overline{P}$ is a violated pattern, and we add it to~$\mathcal{B}$. If $\mathcal{B} \neq \emptyset$ at the end, then the dual solution $(\vpid, \vrhod)$ is infeasible for LM, and we can add some (or all) patterns of $\mathcal{B}$ to RLM and re-optimize it. Otherwise, the dual solution $(\vpid, \vrhod)$ is feasible for LM\@.

An interesting question is how this algorithm can efficiently find a violated pattern without employing dominance rules such as the LSA\@. The answer lies in the ordering of $\overline{I}$. In our instances, there are typically only a few items involved in conflicts or cuts. Additionally, when we are in a state $(i, r)$ with $i \leq |\overline{I_3}|$, $dp(i, r)$ precisely computes the solution value for the current partial pattern. Therefore, although our algorithm may be practically an enumeration process over the items $\overline{I_1}$ and $\overline{I_2}$, these sets are relatively small, causing no performance issues. This allows us to avoid using dominance rules, which could negatively impact pattern diversification. Furthermore, prioritizing processing items belonging to SRIs further enhances efficiency since we observe that SRIs tend to weaken the lower bound of $dp$ more than conflicts do.

The algorithm, as described, identifies all violated patterns. To focus on choosing only a diversified subset, we restrict the exploration of the right branch only if each item $i$ in the partial pattern $\overline{P}$ appears fewer than $2 \zeta$ times in $\mathcal{B}$, with $\zeta$ empirically set to 3. Choosing more than one pattern for each item is beneficial, but a high quantity is redundant. By using this pruning strategy, the ordering of $\overline{I}$ inherently prioritizes patterns with larger items. However, as the algorithm may still traverse many paths, we introduce a stopping criterion: if $\mathcal{B} \neq \emptyset$ and the number of recursive calls $N^R$ exceeds a constant $N^R_{\max}$.
    
In the end, we filter $\mathcal{B}$ by excluding patterns with the highest reduced cost to ensure each item appears at most $\zeta$ times. This strategy compensates for the fact that patterns are not necessarily found in descending order of reduced cost, as we prioritize exploring the left branch over processing the state with the best lower bound. Algorithm~\ref{alg::enumeration} outlines the pseudocode for this approach. We assume that the table $dp$ is precomputed and use global variables $\mathcal{B}$, $dp$, $\overline{I}$, and $N^R$ for simplicity. In our experiments, we set $N^R_{\max} = \frac{|\overline{I}| \cdot W}{10}$ to achieve good-quality results without significantly increasing runtime.

\begin{algorithm}[!t]
    \small
    \SetKwFunction{MultiplePatternGeneration}{MultiplePatternGeneration}
    \SetKwInOut{Input}{input}\SetKwInOut{Output}{output}
    \Input{Partial pattern $\overline{P}$, index $i$, wasted capacity $r$ in $\overline{P}$}

    \If{$i$ = 0 or $r$ = 0}{
        
        Add $\overline{P}$ to $\mathcal{B}$ \hspace{0.5cm} \tcp{ We found a pattern with $\overline{c}_{\overline{P}}(\vpii, \vrhoi, K) < \frac{-K}{M}$}
        
        \textbf{return}
        }
        

    $N^R \gets N^R + 1$
    
    \If{$N^R > N^R_{\max}$ \textbf{and} $\mathcal{B} \neq \emptyset$}{
        
        \textbf{return} \hspace{0.5cm} \tcp{$(\vpid, \vrhod)$ is not feasible, and we reach the iteration limit}
        
        }
        

        \If{$w_{\overline{I}[i]} \leq r$ and $\overline{I}[i]$ does not have a conflict with any item in $\overline{P}$}{ 

            Add $\overline{I}[i]$ to $\overline{P}$ and decrease $r$ by $w_{\overline{I}[i]}$
            
            \If{$dp[i - 1][r] - v(\overline{P}) < \frac{-K}{M}$}{
                
                \MultiplePatternGeneration{$\overline{P}, i - 1, r$}
                }
                Remove $\overline{I}[i]$ from $\overline{P}$ and increase $r$ by $w_{\overline{I}[i]}$
                
                \If{there is an item $j \in \overline{\mathcal{P}}$ that appears $2 \zeta$ or more times in $\mathcal{B}$}{
                    
            \textbf{return} \tcp{Because we want diversity}
            }
        }

        
        $i \gets \max(\{i' : i' < i$ \textbf{and} ${ w_{\overline{I}[i']} \neq w_{\overline{I}[i]} }\})$ \tcp{Try next item with a different size}

    
    \If{$dp[i][r] - v(\overline{P}) < 0$}{
        
        \MultiplePatternGeneration{$\overline{P}, i, r$}
        }
        \caption{\label{alg::enumeration}MultiplePatternGeneration.}
\end{algorithm}
    
We notice that adding the most violated pattern is crucial for fast convergence. Since our algorithm does not always find it, we conduct a second search. This involves maintaining an incumbent pattern~$\overline{P}^\text{inc}$, initialized with the best pattern in $\mathcal{B}$. Thus, we employ a B\&B approach similar to Algorithm~\ref{alg::enumeration}, but explore a branch only if it can produce a better incumbent pattern. We still use the limit of recursive calls $N^R_\text{max}$, but in this case, we prioritize the branch with the smallest lower bound rather than the left branch.

Finally, the alternative pricer mentioned in Section~\ref{sec::num_safe} for determining a lower bound \( B \) for \(\overline{c}(\vpii, \vrhoi, K)\) uses a heap-based approach. It replaces the recursive function's arguments with labels $(\overline{P}, i, r)$ and processes the label with the smallest lower bound (which no longer needs to be smaller than $\frac{-K}{M}$) that can improve \(\overline{P}^\text{inc}\). Note that \( B \) is the minimum of the reduced cost of \(\overline{P}^\text{inc}\) and the lower bound of the unprocessed labels. For performance, we use an iteration limit \( N^R_\text{max} = \frac{|\overline{I}| W}{50} \) and also halt if \( B \geq -K \cdot 10^{-13} \). In this case, the difference between \( z^\text{safe} \) obtained using \( B \) and \(\overline{c}(\vpii, \vrhoi, K)\) is negligible.

\subsection{Generation of Subset-Row Inequalities\label{sec::cutting_generation}}
We generate SRIs by enumeration, which naively involves checking all triples $(i, j, k)$, resulting in $O(n^3)$ triple analysis. Fortunately, it is possible to have a better result in practice. Given a floating-optimal primal solution $\vlambda$ for RLM, let $\overline{I}' = \{i : i \in \overline{I} \wedge d_i = 1\}$ be the set of items with unitary demand, and $\delta_{ij}$ be the affinity between items $i$ and~$j$ defined in Section~\ref{sec::scheme}, with $i, j \in \overline{I}'$.
\begin{lemma}
    Given a triple $S = \{i, j, k\}$ of distinct items, $S$ induces a violated clique cut that can improve the lower bound only if $\delta_{ij} + \delta_{jk} + \delta_{ki} > 1$ and at least two terms of $\delta_{ij} + \delta_{jk} + \delta_{ki}$ are greater than 0\@.
\end{lemma}
\textbf{Proof.} By definition, $S$ induces a violated clique cut if $ \sum_{P \in \mathcal{P}}  \lfloor \frac{a_i^P + a_j^P + a_k^P}{2}  \rfloor \lambda_P^* > 1$. Note that $\delta_{ij} + \delta_{jk} + \delta_{ki} = \sum_{P \in \mathcal{P}} (a_i^P a_j^P + a_i^P a_k^P + a_j^P a_k^P)\lambda_P^* \geq \sum_{P \in \mathcal{P}}  \lfloor \frac{a_i^P + a_j^P + a_k^P}{2}  \rfloor \lambda_P^* > 1$, since, if the coefficient on the right side is $1$, there are at least two items of the triple in $P$, ensuring the left side is at least $1$. Now, suppose $S$ is a cut such that $\delta_{ij} + \delta_{jk} + \delta_{ki} > 1$ and only $\delta_{ij} > 0$, without loss of generality. This implies that $\delta_{ij} > 1$, which can only occur if $i$ and $j$ are over-covered. Since we can convert any (fractional) solution with over-covered items into an exact-covered solution with the same value, $S$ cannot improve the lower bound.~\hfill$\blacksquare$

Let $F_i$ be the adjacency list of each item $i$, where $j \in F_i$ if $\delta_{ij} > 0$. Thus, for each $i \in \overline{I}'$, we analyze all pairs of items $(j, k)$ belonging to $F_i$, checking if $(i, j, k)$ forms a violated cut only when $\delta_{ij} + \delta_{ik} + \delta_{jk} > 1$. If there are many items by pattern, we still need to analyze $O(n^3)$ triples in the worst case. However, this separation algorithm is used intensely only in the hardest instances, where there are few items in each pattern, and, thus, the number of triples analyzed is typically much closer to $O(n^2)$ than $O(n^3)$.

In each node of the B\&B tree, if its safe bound $z^\text{safe}$ allows the existence of an improvement solution, then we trigger the cutting plane generator to seek violated SRIs. If more than $\beta$ violated cuts are identified, we add only the most violated $\beta$ to the LM\@. We limit the generation of cuts to $\alpha$ iterations at each node, stopping regardless of whether any violated cuts remain. We use these limitations because adding each cut makes the model slow to solve, and transitioning to the branching stage is often better than adding cuts for many iterations. In our experiments, the best values for these constants were $\alpha = 10$ and $\beta = 20$.

\section{Primal Heuristics\label{sec::PH}}

For the CSP, we use the Best Fit Decreasing~(BFD) heuristic to produce the initial incumbent solution, as it usually gives floating-optimal solutions.  However, our branching rule focuses solely on enhancing the lower bound, as explained later, making it unlikely to improve the integrality of the relaxation and achieve an integer solution. Thus, we rely on powerful primal heuristics to refine the upper bound, ensuring compatibility with our CRG-based formulation. This compatibility is crucial since RLM often lacks sufficient patterns to produce a better incumbent solution. Hence, effective heuristics should generate new patterns, with any added constraints seamlessly integrated into the pricer. 

Given that the heuristics used in previous works based on SCF~\citep{Belov_2006,Wei_2020,Pessoa_2020} are not competitive with the current state-of-the-art~\citep{Loti_2022}, we decide to use different heuristics, detailed in the following sections. Note that we can transform a (fractional) solution $\vlambda$ of a non-binary model into a binary model solution by dividing each pattern $P$ with $\overline{\lambda}_P > 0$ to $\lfloor \overline{\lambda}_P \rfloor$ patterns with value $\overline{\lambda}'_P = 1$ and one pattern with value $\overline{\lambda}'_P = \overline{\lambda}_P - \lfloor \overline{\lambda}_P \rfloor$. Hereafter, we consider $\vlambda$ as a solution for a binary model in this section.

\subsection{Rounding Heuristic\label{sec::rounding}}
The Rounding heuristic is an inexpensive heuristic that rounds the relaxation and uses the BFD heuristic to pack the possible residual instance. Given the incumbent solution $\Zb$, the total space of an improvement solution $\Zint$ is at most ${(z(\Zb) - 1) \cdot W}$, implying that the total waste in $\Zint$ is less than or equal to the limit ${\Rt = {(z(\Zb) - 1) \cdot W} - \sum_{i \in I} d_i w_i}$. Next, we construct a heuristic solution $S$ by using the current RLM solution $\vlambda$, a variable $\Rh$ indicating the residue that $S$ can have to be an improvement solution (where $\Rh = \Rt$ at the beginning) and a threshold $\lambda_{\min} = 0.6$.

We build $S$ by sequentially iterating over patterns in a non-increasing order of value. For each pattern $P$, we check if $\overline{\lambda}'_P \geq \lambda_{\min}$ and if $\Rh \geq R_P$. If both conditions are met, we add $P$ to $S$ and update $\Rh \leftarrow \Rh - R_P$. Notably, if utilizing pattern $P$ results in over-covering certain items, the excess demand is considered as waste. After evaluating all patterns, there may be a subset $\overline{I}'$ of non-packed items, which we pack using the BFD heuristic. Given its low computational cost, we apply this heuristic to every feasible RLM solution $\lambda'$, encompassing intermediate solutions of the CRG process.

\subsection{Relax-and-Fix Heuristic}
The Relax-and-Fix~(RF) heuristic is a diving heuristic~\citep{Belvaux_2000} that seeks an integer solution through an interactive process of relaxing and fixing variables. At each step, it optimally solves the relaxed model and selects a subset of patterns $F \subseteq \mathcal{P}$, fixing their values to an integer value. Using a strong safe lower bound $z^\text{safe}$ in the CSP as that provided by the SCF, it is often easy to find an incumbent solution $\Zb$ such that ${z(\Zb) - \lceil z^\text{safe} \rceil = 1}$. However, closing this gap can prove challenging. Given the tightness of this gap, we get better results by fixing a small set of variables $F$ in each step.

At each step, our algorithm optimizes RLM using only a subset of items $\overline{I}'$, initially set as $\overline{I}$, and let $\vlambda$ be the primal solution obtained when the CRG process halts. First, if there are variables $P$ with $\vlambda_P = 1$, we choose $F$ as these variables. Otherwise, we calculate the gap $g = z(\Zb) - z^\text{safe} - 1$ and choose $F$ by iterating over the variables $P \in \vlambda$ in a non-increasing order of value $\vlambda_P$. The first pattern is added to $F$ to prevent looping, and each subsequent pattern $P$ is added to $F$ if $\vlambda_P > 0.5$, $(1 - \vlambda_P) \leq g$, and the set of over-covered items remains unchanged. We update $g$ whenever a pattern $P$ is added to $F$ by subtracting $1 - \vlambda_P$. Finally, upon selecting $F$, we add the patterns of $F$ to a partial solution $S$ and remove their corresponding items from $\overline{I}'$. This process repeats until $\overline{I}'$ becomes empty, indicating that $S$ is a feasible solution.

For performance reasons, we decided not generate new cuts in this heuristic and to halt the CG process when the current $\vlambda$ is sufficient to yield a better $\vlambda^\text{inc}$, that is, when $z(\vlambda) + |F| \leq z(\vlambda^\text{inc}) - 1$. The Rounding heuristic is executed after each CG iteration aiming to find an integer solution for $\overline{I}'$. Furthermore, we execute the algorithm three times, incorporating a diversification strategy. In the second and third runs, we start with a subset of the prior fixed solution $S$, constructed by excluding the last quarter of the selected sets $F$. This approach aims to mitigate potentially suboptimal choices in $S$, which are more likely to occur towards the end. Finally, we trigger this heuristic every ten left branches, starting from the root node.

\subsection{Constrained Relax-and-Fix Heuristic}
Given a floating-optimal (vertex) solution $\vlambda$ for the root node without SRIs, we define the \emph{Polyhedron Integrality Ratio} of LM as $\sum_{P \in \mathcal{P}: \vlambda_P \in \mathbb{Z}} \vlambda_P / z(\vlambda)$. We define it considering a single solution (the solution found by the FLPS) because we observe that different solutions under these conditions generally produce similar values. We denote a polyhedron as \emph{highly fractional} if its integrality is smaller than 30\%. Heuristics such as RF, which start with a linear relaxation solution, may struggle in these cases due to the high likelihood of setting incorrect values for variables. To address this, some heuristics leverage an incumbent solution to confine the range of some variables. Local Branching~(LB) proposed by~\citet{Fischetti_2003}, is one such heuristic. Given an incumbent solution $\vlambda^\text{inc}$ and a small integer~$k$, LB solves the resulting ILP after adding the following constraint:
\begin{equation}
\sum_{P \in \mathcal{P}: \vlambda_P^\text{inc} = 0} \lambda_P + \sum_{P \in \mathcal{P}: \vlambda_P^\text{inc} = 1} (1 - \lambda_P) \leq k.
\label{eq::local_branching1}
\end{equation}

However, an improvement solution and the incumbent solution can diverge in many patterns in the CSP\@. On the other hand, perhaps the RF heuristic fails to yield a better solution because it fixes a few incorrect patterns. To address this, we replace the aforementioned constraint with the following:
\begin{equation}
\sum_{P \in S^\text{inc}} \lambda_P \geq  \vert S^\text{inc} \vert  -  k,
\label{eq::local_branching2}
\end{equation}
where $S^\text{inc}$ is the largest set of fixed variables among all runs of the RF such that the resulting lower bound is lower than the incumbent solution value. When all patterns in $S^\text{inc}$ are added to RLM, the pricing problem remains unchanged. However, employing this inequality in a B\&B approach can be inefficient, particularly when dealing with a highly fractional polyhedron, where branching becomes almost ineffective. Therefore, we propose the \emph{Constrained Relax-and-Fix}~(CRF) heuristic, which involves running the RF heuristic with Constraint (\ref{eq::local_branching2}) but without any branch constraints. This means temporarily removing all conflicts and unmerging all items before executing the heuristic, while keeping any still valid cut.

This heuristic runs every 30 left branches if the tree's height is lower than 30, and every 20 left branches otherwise. Moreover, it is enabled to run only after two RF runs, as it depends on the $S^\text{inc}$ generated by RF\@. The frequency intensification is because a high height probably indicates difficulty in improving the incumbent solution. If ten consecutive attempts fail to improve using $S^\text{inc}$, it is discarded in favor of the next best incumbent to prevent being stuck in an unproductive search space. Finally, we perform two executions each time using $k = 6$ and $k = 12$, with each execution performing three RF runs as described earlier.

\section{Model's Optimizations\label{sec::Optimization}}
This section introduces optimizations to enhance our framework's performance. We begin by discussing our branching rule and the technique \emph{splay operation}. Next, we outline our dual inequalities, along with strategies for binary pricing, waste optimization, and model cleaning based on reduced cost.

\subsection{Branching Rule\label{sec::historic_rule}}

In the CSP, we could expand upon the rule used by~\citet{Wei_2020} for the BPP, branching on pair of items $(i, j)$ where $\delta_{ij} \notin \mathbb{Z}$ and is furthest from the nearest integer value (with tie-breaking based on maximizing $w_i + w_j$). However, we can enhance this approach by leveraging historic branching data to prioritize item pairs, a strategy uncommon in branching strategies for cutting and stock problems. In this improved strategy, we choose to branch on the pair of items $(i, j)$ that maximizes the lexicographical order of $(-r_{ij}, w_i + w_j)$, with $\delta_{ij} \notin \mathbb{Z}$, where $r_{ij}$ represents the priority of the pair of items $(i, j)$, which is defined next.

As we notice that a gradual adjustment in the priority produces better results, we employ a list-based approach. Let $s$ denote a sequence of item pairs, initially empty. We define $r_{ij}$ as the index of $(i, j)$ in $s$ if $(i, j)$ appear in $s$, and $r_{ij} = \infty$, otherwise. We update $s$ to prioritize pairs of items that likely cannot belong to the same or different patterns, making them promising candidates to increase the lower bound. Suppose that the pair of items $(i, j)$ was chosen to branch at the current B\&B node. In the case where both branches are pruned by bound, we append $(i, j)$ to the end of $s$, if $(i, j) \notin s$, and we swap $(i, j)$ with the previous element of $s$~(if it exists), otherwise. And in the case where the left branch cannot be pruned by bound and $(i, j) \in s$,  we swap $(i, j)$ with the next element of $s$ (if it exists). In any other case, $s$ remains unchanged.

Note that pairs of items chosen for branching at deeper levels of the B\&B tree may frequently not be good choices near the root node. This is because these pairs might have been chosen simply because the ancestor branches of these items, which remained fixed during iterations, did not allow for an improvement solution. Our list-based approach addresses this by making the priority of these pairs of items increase gradually.

\subsection{Splay Operation}

We explore the B\&B tree by performing a depth-first search~(DFS) and propose a strategy called \emph{splay operation} to control its height by discarding left branch nodes and avoiding deep dives that generally do not contribute to increasing the lower bound. 

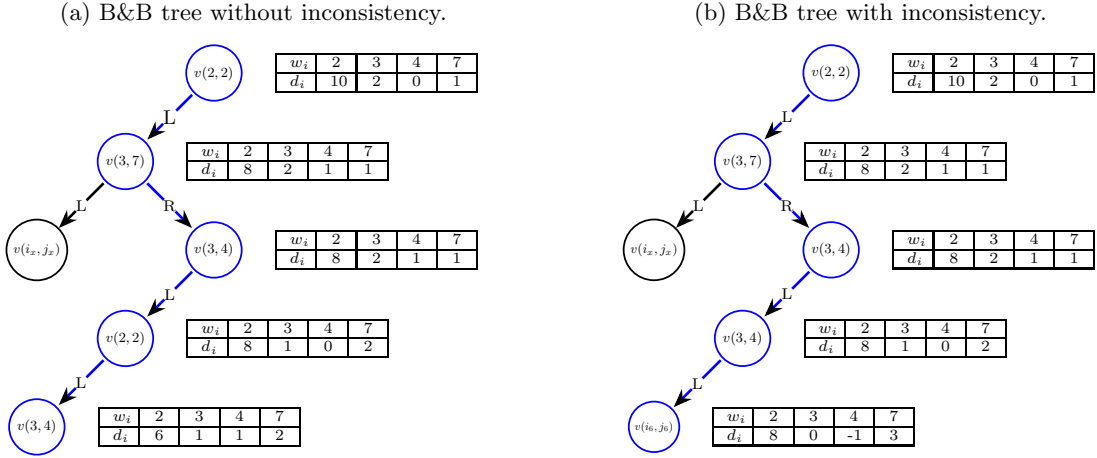
\begin{figure}[!t]
    \centering
    \caption{Examples of trees with and without inconsistency. The blue color indicates the active path in the B\&B tree.}
    \begin{subfigure}{.49\textwidth}
        \caption{B\&B tree without inconsistency.}
        \label{fig::tree}
        \centering
        \resizebox{0.8\textwidth}{!}{%
\begin{tikzpicture}
\begin{scope}[every node/.style={circle,thick,draw = blue, outer sep=2pt}, every label/.style={shape=rectangle, draw=none, fill=none}]
\node [label = right:{ \hspace{0.3cm}\scriptsize
    \begin{tabular}{|C{0.2cm}|C{0.2cm}|C{0.2cm}|C{0.2cm}|C{0.2cm}|}
        \hline
        $w_i$ & 2  & 3 & 4 & 7 \\ \hline
        $d_i$ & 10 & 2 & 0 & 1 \\ \hline
    \end{tabular}
}, scale=0.7](A) at (0,0) { \small $v(2, 2)$};

\node [label = right:{ \hspace{0.3cm}\scriptsize
    \begin{tabular}{|C{0.2cm}|C{0.2cm}|C{0.2cm}|C{0.2cm}|C{0.2cm}|}
        \hline
        $w_i$ & 2  & 3 & 4 & 7 \\ \hline
        $d_i$ & 8 & 2 & 1 & 1 \\ \hline
    \end{tabular}
}, scale=0.7](B) at (-1.4,-1.4) { \small $v(3, 7)$};

\node [label = right:{ \hspace{0.3cm}\scriptsize
    \begin{tabular}{|C{0.2cm}|C{0.2cm}|C{0.2cm}|C{0.2cm}|C{0.2cm}|}
        \hline
        $w_i$ & 2  & 3 & 4 & 7 \\ \hline
        $d_i$ & 8 & 2 & 1 & 1 \\ \hline
    \end{tabular}
}, scale=0.7] (C) at (0,-2.8) { \small $v(3, 4)$};

\node [label = right:{ \hspace{0.3cm}\scriptsize
    \begin{tabular}{|C{0.2cm}|C{0.2cm}|C{0.2cm}|C{0.2cm}|C{0.2cm}|}
        \hline
        $w_i$ & 2  & 3 & 4 & 7 \\ \hline
        $d_i$ & 8 & 1 & 0 & 2 \\ \hline
    \end{tabular}
}, scale=0.7] (D) at (-1.4,-4.2) { \small $v(2, 2)$};

\node [label = right:{ \hspace{0.3cm}\scriptsize
    \begin{tabular}{|C{0.2cm}|C{0.2cm}|C{0.2cm}|C{0.2cm}|C{0.2cm}|}
        \hline
        $w_i$ & 2  & 3 & 4 & 7 \\ \hline
        $d_i$ & 6 & 1 & 1 & 2 \\ \hline
    \end{tabular}
}, scale=0.7] (E) at (-2.8,-5.6) { \small $v(3, 4)$};
\end{scope}

\begin{scope}[every node/.style={circle,thick,draw = black}]
    \node[scale=0.65] (F) at (-2.8,-2.8) { \small $v(i_x, j_x)$};
\end{scope}

\begin{scope}[>={Stealth[black]},
              every node/.style={fill=white,circle,inner sep=0pt, font=\scriptsize},
              every edge/.style={draw=blue,very thick}]
    \path [->] (A) edge node {\small L} (B);
    \path [->] (B) edge node {R} (C);
    \path [->] (C) edge node {L} (D);
    \path [->] (D) edge node {L} (E);
\end{scope}

\begin{scope}[>={Stealth[black]},
              every node/.style={fill=white,circle, inner sep=0pt, font=\scriptsize},
              every edge/.style={draw=black,very thick}]
    \path [->] (B) edge node {L} (F);
\end{scope}
\end{tikzpicture}
}
    \end{subfigure}
    \begin{subfigure}{.49\textwidth}
        \caption{B\&B tree with inconsistency.}
        \label{fig::tree2}
        \centering
        \resizebox{0.8\textwidth}{!}{%
\begin{tikzpicture}
\begin{scope}[every node/.style={circle,thick,draw = blue, outer sep=2pt}, every label/.style={shape=rectangle, draw=none, fill=none}]

\node [label = right:{ \hspace{0.3cm}\scriptsize
    \begin{tabular}{|C{0.2cm}|C{0.2cm}|C{0.2cm}|C{0.2cm}|C{0.2cm}|}
        \hline
        $w_i$ & 2  & 3 & 4 & 7 \\ \hline
        $d_i$ & 10 & 2 & 0 & 1 \\ \hline
    \end{tabular}
}, scale=0.7](A) at (0,0) { \small $v(2, 2)$};

\node [label = right:{ \hspace{0.3cm}\scriptsize
    \begin{tabular}{|C{0.2cm}|C{0.2cm}|C{0.2cm}|C{0.2cm}|C{0.2cm}|}
        \hline
        $w_i$ & 2  & 3 & 4 & 7 \\ \hline
        $d_i$ & 8 & 2 & 1 & 1 \\ \hline
    \end{tabular}
}, scale=0.7](B) at (-1.4,-1.4) { \small $v(3, 7)$};

\node [label = right:{ \hspace{0.3cm}\scriptsize
    \begin{tabular}{|C{0.2cm}|C{0.2cm}|C{0.2cm}|C{0.2cm}|C{0.2cm}|}
        \hline
        $w_i$ & 2  & 3 & 4 & 7 \\ \hline
        $d_i$ & 8 & 2 & 1 & 1 \\ \hline
    \end{tabular}
}, scale=0.7] (C) at (0,-2.8) { \small $v(3, 4)$};

\node [label = right:{ \hspace{0.3cm}\scriptsize
    \begin{tabular}{|C{0.2cm}|C{0.2cm}|C{0.2cm}|C{0.2cm}|C{0.2cm}|}
        \hline
        $w_i$ & 2  & 3 & 4 & 7 \\ \hline
        $d_i$ & 8 & 1 & 0 & 2 \\ \hline
    \end{tabular}
}, scale=0.7] (D) at (-1.4,-4.2) { \small $v(3, 4)$};

\node [scale=0.85, label = right:{ \hspace{0.3cm}\scriptsize
    \begin{tabular}{|C{0.2cm}|C{0.2cm}|C{0.2cm}|C{0.2cm}|C{0.2cm}|}
        \hline
        $w_i$ & 2  & 3 & 4 & 7 \\ \hline
        $d_i$ & 8 & 0 & \mbox{-1} & 3 \\ \hline
    \end{tabular}
}, scale=0.65] (E) at (-2.8,-5.6) { \small $v(i_6, j_6)$};
\end{scope}

\begin{scope}[every node/.style={circle,thick,draw = black}]
    \node[scale=0.65] (F) at (-2.8,-2.8) { \small $v(i_x, j_x)$};
\end{scope}

\begin{scope}[>={Stealth[black]},
              every node/.style={fill=white,circle,inner sep=0pt, font=\scriptsize},
              every edge/.style={draw=blue,very thick}]
    \path [->] (A) edge node {L} (B);
    \path [->] (B) edge node {R} (C);
    \path [->] (C) edge node {L} (D);
    \path [->] (D) edge node {L} (E);
\end{scope}

\begin{scope}[>={Stealth[black]},
              every node/.style={fill=white,circle,inner sep=0pt, font=\scriptsize},
              every edge/.style={draw=black,very thick}]
    \path [->] (B) edge node {L} (F);
\end{scope}
\end{tikzpicture}
}
    \end{subfigure}
\end{figure}

Let $C = (v(i_1, j_1), O_1, v(i_2, j_2), O_2, \dots, O_{t-1}, v(i_t, j_t))$ be the active path in the B\&B tree, where $v(i, j)$ represent a node of branching on the pair of items $(i, j)$, and $O_k \in \{L, R\}$, for ${1 \leq k < t}$, is an edge, where $L$ represents a left branch and $R$ a right branch. Figure~\ref{fig::tree} shows an example with ${C = (v(2, 2), L, v(3, 7), R, v(3, 4), L, v(2, 2), L, v(3, 4))}$, where the item size is the item identifier, and the tables show the demand list at each node.

We say that a contiguous sequence $(v(i_k, j_k), v(i_{k + 1}, j_{k + 1}), \ldots, v(i_l, i_l))$, with $1 \leq k \leq l \leq t$ is a prefix (suffix) of vertices from $C$ if $k = 1$ ($l = t$). Let ${\mathcal{V}_\text{suff} = (v_1, \ldots, v_{ \vert \mathcal{V}_\text{suff} \vert })}$ be the maximal suffix of vertices connected only by $L$ edges, where ${\mathcal{V}_\text{suff} = (v(i_3, j_3), v(i_4, j_4), v(i_5, j_5))}$ in the example. A priori, we can discard any subset of $\mathcal{V}_\text{suff}$ without losing the algorithm's correctness since we follow the DFS order. The only exception is when we remove an intermediary node that creates an item used by a posterior node. In Figure~\ref{fig::tree2}, see that after removing node $v(i_4 = 2, j_4 = 2)$, the resulting tree is \emph{inconsistent}, since the node $v(i_5 = 3, j_5 = 4)$ has demand zero in the item of size $4$ and it produces a negative demand in its left child $v(i_6, j_6)$.

At each B\&B node $v(i, j)$ where both child nodes' relaxations fail to improve the incumbent solution, we execute the splay operation, detailed below. The splay operation involves moving $v(i, j)$ closer to the root node by removing a sequence of $L$ edges. To determine the nodes to be removed, we construct a list $V_r$ by iterating over $k$ from $|\mathcal{V}_{\text{suff}}| - 1$ down to $1$. We add $v_k$ to $V_r$ only if it is not fixed and maintains tree consistency, meaning that removing the nodes in ${v_k} \cup V_r$ results in all prefixes of the resulting active path $C$ having non-negative demand for all items. After that, if $V_r \neq \emptyset$, we remove its nodes, mark $v(i, j)$ as non-removable (to prevent cyclic processing in following iterations), and reprocess this node. Otherwise, the splay operation fails, and we proceed to the next node in the DFS order.

Finally, in the DFS order prioritizing the left branch, exploring the right branch of a node $v(i, j)$ implies that there is no improvement solution where $i$ and $j$ belong to the same patterns while maintaining the decisions of $v(i, j)$'s ancestors, as such a solution should have been found when exploring the left branch. Conversely, exploring the left branch does not provide any definitive implication about $i$ and $j$. Since the splay operation only removes nodes from $\mathcal{V}_\text{suff}$, it is not a pruning strategy. However, it helps raise the lower bound by moving items that cannot belong to the same or different patterns closer to the root node.

\subsection{Dual Inequalities and Binary Pricing}
The Dual Inequalities and Binary Pricing techniques discussed below aim to solve the root node without any SRIs, thus we omit $\vrho$ in dual solutions. Dual inequalities \citep{Amor_2006, Irnich_2016} offer another promising strategy to address the slow convergence to solve the root node in CG-based formulations. These studies reveal that there exists a dual optimal solution $\vpi^*$ for LM, where the value $\vpi_i^*$ of each item $i$ is proportional to its size $w_i$. Consequently, stabilizing the dual values $\vpi$ can be done by adding constraints of the form ${\vpi_i \geq \sum_{j \in S} t_j \vpi_j}$ where ${w_i \geq \sum_{j \in S} t_s w_j}$, $S \subseteq \overline{I} \setminus{i}$, and $t \in \mathbb{Z}^S_{> 0}$. Given the superpolynomial number of such inequalities, these authors often employ fixed-size subsets or dynamically separate them. Importantly, these inequalities do not compromise relaxation strength, and in cases where they yield an infeasible primal solution, conversion to a feasible one is straightforward. 

Our approach deviates from this by introducing dual constraints that deliberately weaken the lower bound of the relaxation. As slow convergence is intuitively justified by the high instability of item values among the CG iterations, our approach focuses on stabilizing them. Given a floating-optimal dual solution $\vpi$ for RLM (which may not be floating-optimal for LM), our dual inequality set is defined as%
\begin{equation}
D = \{\pi_i \leq \gamma w_i : \forall i \in \overline{I}\},\label{eq::dual_inequalities}
\end{equation}
where $\gamma$ is the value of a unit of size in solution $\vpi$, i.e.,  $\gamma = \frac{z(\vpi)}{\sum_{i \in \overline{I}} d_i w_i}$. These inequalities constrain the cost-benefit $\pi_i$ of each item $i$ to be at most $\gamma$. However, each item $i$ corresponds to a variable from \eqref{eq::dual_inequalities} in the primal model, which covers $i$ with an objective cost of $\gamma w_i$. These variables are distinct from patterns, as patterns have a cost of 1, whereas $\gamma w_i$ is typically less than 1. Consequently, \eqref{eq::dual_inequalities} are not valid dual inequalities for the original problem, potentially resulting in $\vpi$ having a lower value than the floating-optimal solution of LM\@. Additionally, we choose to update $D$ only when there is no item~$i$ with $\vpi_i = \gamma w_i$, i.e., when all inequalities $D$ are strict. Alternatively, we could update it whenever any inequality becomes strict, but we have observed that this often results in a very weak lower bound.

At the end of CG, the solution found may be floating-optimal only for the relaxed problem with $D$, not for the original relaxation. So, we remove $D$ and generate more columns for optimality. This re-optimization was typically fast in our experiments, often requiring only a few additional columns.

Finally, we can speed up this first root optimization using binary pricing, as used by~\citet{Delorme_2020}, which limits the pricer to generate patterns with only one item of each distinct size. Despite this restriction, dual stabilization still performs well, with most iterations occurring in the first stage, where items have unitary demand, making it quicker. Then, the second root optimization, which occurs without the dual inequalities and the binary pricing, requires only a few additional columns to converge. Notably, the binary pricing strategy is beneficial primarily in instances with high item demands. Based on our experiments, we decided to enable this strategy only in instances where $\vert I \vert < \frac{5}{6} \sum_{i \in I} d'_i$.

\subsection{Waste Optimization}
Given the total waste $\Rt$ defined in Section~\ref{sec::rounding}, it is important to note that every pattern in an improvement solution exhibits waste $R_P \leq \Rt$. Hence, during the B\&B process, patterns violating this limit can be excluded, and the pricer can be modified to generate patterns with waste at most $\Rt$. This modification involves adjusting the DP basis case: setting $dp(0, r) = K$ if $r \leq \Rt$ and $dp(0, r) = \infty$ otherwise. This approach, used by \citet{Loti_2022,Delorme_2020} after solving the root node, is extended by us in Appendix~\ref{app::waste} to apply while solving the root node, reducing solution time for the most costly node.

\subsection{Model Cleaning by Reduced Cost}
Reduced Cost Variable Fixing~(RCVF) is a technique widely utilized in literature~(see, e.g.,~\cite{Crowder_1983}), which achieves high effectiveness mainly in arc-flow models~(\citet{Delorme_2020, Loti_2022}). Let $\vpi$ be a feasible dual solution, $\overline{\lambda}_{\mathbb{Z}}$ be a feasible integer solution, and $\overline{c}_P$ be the reduced cost of a pattern $P$ in $\vpi$. By the reduced cost definition, a pattern $P$ such that ${z(\vpi) + \overline{c}_P > z(\overline{\lambda}_{\mathbb{Z}}) -  1}$ cannot belong to a primal integer solution better than $\overline{\lambda}_{\mathbb{Z}}$. Thus, we can use this strategy to fix variables to $0$. 

In arc-flow models, this idea is highly effective due to their typically pseudo-polynomial number of arcs. In CG-based models, arcs can be directly fixed to $0$ in the pricing problem without altering its structure. However, Set Covering models lack these advantages, often having an exponential number of variables and requiring a list of forbidden patterns in the pricer. Moreover, RCVF is typically applied only at the root node in the literature. Nonetheless, we can extend this technique to other B\&B nodes, keeping the variables fixed only within the node's subtree.

Given the high cardinality of $\mathcal{P}$ and the relatively low number of fixed patterns, our approach removes patterns from RLM that cannot improve the incumbent solution but allow their potential reintroduction later.  This method is called Model-Cleaning-by-Reduced-Cost (MCRC). Both MCRC and Waste Optimization are performed as a preprocessing step before optimizing any B\&B node.



\section{Computational Experiments for the CSP\label{sec::CE}}
Next, we present our computational experiments\footnote{All instances and the source code of our framework for CSP and the other problems are available at \href{https://gitlab.com/renanfernandofranco/a-branch-and-cut-and-price-algorithm-for-cutting-stock-and-related-problems}{ https://gitlab.com/renanfernandofranco/a-branch-and-cut-and-price-algorithm-for-cutting-stock-and-related-problems}.}
for the CSP, comparing our framework to other state-of-the-art algorithms and showing the results for our benchmark. Since some challenging instances for BPP also are challenging for CSP, this comparison also involves BPP's algorithms. Due to space constraints,  tests demonstrating the usefulness of our algorithm's features are presented in Appendices~\ref{app::features_tests} and~\ref{app::exp1h}. In this article, all results for our framework were obtained with a computer running operational system Ubuntu 18.04.6 LTS (64 bits), using the language C++17 with the compiler GCC 7.5.0, and processor Intel\textsuperscript{\textregistered}\ Xeon\textsuperscript{\textregistered}\ CPU E5-2630 v4 @ 2.20GHz with 64 GB of RAM, which has a single-thread passmark indicators (STPI) equal to $1785$. We employ these indicators in subsequent sections to compare CPU performance, available at \href{www.passmark.com}{www.passmark.com}, where higher values denote better performance.

Moreover, we use the floating-point solver Gurobi 10.0.3 executed in a \textit{single thread} to solve the RLM, usually using the simplex method since it uses the previous basis as a warm start for the current optimization. Exceptions arise when addressing the numerical issue outlined in Section~\ref{sec::num_safe} and after adding numerous patterns to the RLM following the execution of the RF heuristic, where we resort to the barrier method.

The following tables have columns \emph{Opt}, \emph{Time}, \emph{Cols}, and \emph{Cuts} that represent the number of instances solved to optimality, the average time in seconds, the average number of columns generated, and the average number of cuts generated, respectively, for each class of instances and algorithms studied.

\subsection{Comparison with state-of-the-art algorithms}
In Table~\ref{Table1}, we compare our algorithm with state-of-the-art algorithms for the BPP/CSP using a time limit of 1 hour: \textsc{Belov} from~(\citet{Belov_2006}), NF-F from~\cite{Loti_2022} and BCCF from~\citet{Baldacci_2024}. The results for 
\textsc{Belov} were obtained with an Intel Xeon at 3.1-GHz (STPI $1544$)\footnote{The algorithm was executed by~\citet{Delorme_2020}.} and 8 GB RAM\@; 
for NF-F were obtained with an Intel Xeon E3-1245 v5 at 3.50GHz (STPI $2249$) and 32 GB RAM\@; 
and for BCCF were obtained with an Intel Xeon Gold 6130 at 2.10 GHz (STPI $2067$) and 192 GB of RAM\@. The instances used for the comparison belong to BPP Lib (\citet{Delorme_2018}), a set of instances proposed in the last decades and used by most of the algorithms in the literature.

Our algorithm surpasses all other state-of-the-art algorithms, solving significantly more instances in notably less time. Specifically, it solves 14 additional AI and ANI instances compared to \textsc{NF-F}, the current state-of-the-art, leaving only 10 AI instances and 3 ANI instances unsolved.

Our algorithm performance is slightly worse than \textsc{Belov} in the Waescher and GI classes. As we can see, we encountered difficulty finding the optimal solution for one instance in the GI-BA class. \textsc{Belov}'s superior performance in these classes can be attributed to its pricing problem, based on fractional knapsack, which is less impacted by large bin capacities, such as those in the GI class where ${W \in \{5 \cdot 10^5, 1.5 \cdot 10^6\}}$. Conversely, \textsc{Belov} exhibits a worse performance in instance classes such as AI and ANI, which require identifying numerous wasteless patterns.

Other algorithms in the literature~(see, e.g., \citet{Delorme_2020,Wei_2020}), omitted here due to space constraints, also effectively addressed all classes in the literature except class AI and ANI\@. However, these two classes pose challenges, with many instances remaining unsolved until the recent publication of NF-F (\citet{Loti_2022}), which successfully tackled instances with 600 items or more. NF-F's efficiency in these classes heavily relies on its diversified generation of multiple patterns, yet this alone proves insufficient as our algorithm requires the smaller tolerance $\epsilon'$ defined in Section~\ref{sec::num_safe} to handle such instances. Thus, we speculate that NF-F also mitigates numerical issues in these large instances, possibly through its arc fixing strategy, which by removing arcs from the pricing provides a reduction in the range of reduced costs.

\begin{table}[H]
\centering
\caption{Comparison with state-of-the-art algorithms for the BPP/CSP (time limit of 3600s).}
\footnotesize
\begin{tabular}{ccccccccccccccccccccc}
\hline
                &                   &  \multicolumn{2}{c}{\textbf{Belov}}   && \multicolumn{2}{c}{\textbf{NF-F}}&& \multicolumn{2}{c}{\textbf{BCCF}} &&  \multicolumn{4}{c}{\textbf{Our}}                              \\ \cline{3-4} \cline{6-7} \cline{9-10} \cline{12-15} 
\textbf{Class}  & \textbf{Total}    &  \textbf{Opt}  & \textbf{Time}        && \textbf{Opt}   & \textbf{Time}   &&  \textbf{Opt}& \textbf{Time} && \textbf{Opt}     & \textbf{Time}     & \textbf{Cols} & \textbf{Cuts} \\
\midrule
AI 202          &    50             &    \textbf{50} &           90.6       &&    \textbf{50} &     2.0         &&          –       &       –        &&    \textbf{50}   &    \textbf{0.4}   &  1000.8       &   10.9 \\
AI 403          &    50             &             45 &          699.4       &&    \textbf{50} &    25.2         &&          –       &       –        &&    \textbf{50}   &    \textbf{5.0}   &  2619.6       &   23.7 \\
AI 601          &    50             &              – &              –       &&             49 &   192.4         &&          –       &       –        &&    \textbf{50}   &   \textbf{57.1}   &  4778.8       &   56.3 \\
AI 802          &    50             &              – &              –       &&             46 &   566.5         &&          –       &       –        &&    \textbf{48}   &  \textbf{223.7}   &  7074.5       &   97.5 \\
AI 1003         &    50             &              – &              –       &&             36 &  1577.1         &&          –       &       –        &&    \textbf{42}   &  \textbf{794.9}   &  9997.0       &  134.6 \\
ANI 201         &    50             &    \textbf{50} &          144.2       &&    \textbf{50} &     3.0         && \textbf{50}      &    13.6        &&    \textbf{50}   &    \textbf{0.4}   &  1034.5       &    8.0 \\
ANI 402         &    50             &              1 &         3555.6       &&    \textbf{50} &    24.9         && \textbf{50}      &   308.2        &&    \textbf{50}   &    \textbf{2.2}   &  2438.8       &    8.2 \\
ANI 600         &    50             &              – &              –       &&    \textbf{50} &   140.7         &&          25      &  1931.5        &&    \textbf{50}   &   \textbf{11.8}   &  4129.5       &   18.7 \\
ANI 801         &    50             &              – &              –       &&             49 &   393.2         &&           3      &  3352.7        &&    \textbf{50}   &   \textbf{57.0}   &  6154.7       &   31.7 \\
ANI 1002        &    50             &              – &              –       &&             43 &  1302.5         &&           0      &  3600.0        &&    \textbf{47}   &  \textbf{374.0}   &  8819.7       &   56.0 \\
FalkenauerT     &    80             &    \textbf{80} &           56.9       &&    \textbf{80} &     0.3         &&           –      &       –        &&    \textbf{80}   &   \textbf{0.1}    &   517.0       &    1.2 \\
FalkenauerU     &    80             &    \textbf{80} &        $< 0.1$       &&    \textbf{80} &     0.1         &&           –      &       –        &&    \textbf{80}   &            0.01   &    88.6       &    0.0 \\
Hard            &    28             &    \textbf{28} &            7.5       &&    \textbf{28} &    23.6         &&           –      &       –        &&    \textbf{28}   &    \textbf{6.4}    &   930.1       &   36.1 \\
GI AA           &    60             &    \textbf{60} &   \textbf{2.8}       &&              – &       –         &&           –      &       –        &&    \textbf{60}   &            14.5   &  1279.7       &    0.0 \\
GI AB           &    60             &    \textbf{60} &           10.9       &&              – &       –         &&           –      &       –        &&    \textbf{60}   &   \textbf{7.9}    &  1016.8       &    0.0 \\
GI BA           &    60             &    \textbf{60} &   \textbf{2.8}       &&              – &       –         &&           –      &       –        &&             59   &           102.0   &  1265.4       &    0.0 \\
GI BB           &    60             &    \textbf{60} &  \textbf{10.5}       &&              – &       –         &&           –      &       –        &&    \textbf{60}   &            25.7   &  1091.0       &    0.0 \\
Random          &  3840             &              – &              –       &&  \textbf{3840} &     0.9         &&           –      &       –        &&  \textbf{3840}   &   \textbf{0.05}   &   213.8       &    0.05 \\
Scholl          &  1210             &  \textbf{1210} &            0.2       &&  \textbf{1210} &     1.4         &&           –      &       –        &&  \textbf{1210}   &   \textbf{0.05}   &   111.3       &    1.2 \\
Schwerin        &   200             &   \textbf{200} &            1.1       &&   \textbf{200} &     0.2         &&           –      &       –        &&   \textbf{200}   &   \textbf{0.02}   &   88.2        &    0.04  \\
Waescher        &    17             &    \textbf{17} &   \textbf{0.1}       &&    \textbf{17} &   161.2         &&           –      &       –        &&    \textbf{17}   &             0.3   &   385.1       &    1.7  \\
\bottomrule
\end{tabular}
\label{Table1}
\end{table}

The BCCF algorithm proposed by~\citet{Baldacci_2024} and introduced after NF-F, addresses a critical concern: the lack of numerical safety in existing algorithms,  particularly due to the reliance on floating-point numbers in pruning phases. This issue is particularly acute in the ANI class, which exclusively comprises non-IRUP instances. To mitigate this, they propose a numerically safe algorithm tailored for the ANI class, which outperforms all existing state-of-the-art algorithms except NF-F. Notably, our algorithm is also numerically safe and surpasses both the NF-F and BCCF algorithms.


\subsection{Our benchmark}

We conjecture that our algorithm easily solves AI instances because their relaxations have a high polyhedron integrality ratio (greater than 50\% in most cases). This indicates that the polyhedron of the linear relaxation of SCF closely resembles the convex hull of integer solutions in these instances. Based on this observation, we propose a benchmark of 50 instances for each ${(n, W) \in {(216, 10^3), (405, 1.5 \cdot 10^3), (648, 2 \cdot 10^3)}}$, with this proportion averaging less than 5\%.

Table~\ref{Table5b} presents the results obtained by our framework in this benchmark. Due to space constraints, we present how to build this benchmark and a more detailed analysis in Appendix~\ref{app::benchmark}. Notably, even with a one-hour time limit, our framework only yielded optimal solutions for 44 out of 50 instances with 216 items, leaving many instances unresolved with 405 and 648 items. Thus, these instances constitute a new and challenging CSP benchmark to be considered in future works and potentially adaptable to other problems.

\begin{table}[ht!]
\centering
\caption{Results for our CSP benchmark using our framework (time limit of 3600s).}
\label{Table5b}
\footnotesize
\begin{tabular}{cccccccccccccc}
\toprule
N  & Total & Opt &  Time &  \mline{Pricing}{Time} & \mline{LP}{Time}  & Cols    &    Cuts\\
\midrule
216 & 50 & 44 & 1082.4 & 220.2 & 730.6 & 3251.3 & 554.6 \\
405 & 50 & 34 & 2026.2 & 552.0 & 1225.9 & 12126.7 & 897.1 \\
648 & 50 & 13 & 3207.3 & 925.5 & 1893.4 & 20078.6 & 1273.6 \\
\bottomrule
\end{tabular}
\end{table}

\section{Application in Other Problems\label{sec::Problems}}

To show the versatility of our framework, we apply it to the following four problems: the Skiving Stock Problem, Identical Parallel Machines Scheduling with Minimum Makespan, the Ordered Open-End Bin Packing Problem, and the Class-Constrained Bin Packing Problem.
We have chosen these problems because they are similar to CSP, although we believe our framework can also be extended to other problems with strong relaxations and fast pricers. Next, we explain the necessary adaptation for these problems and computational experiments between our framework and state-of-the-art algorithms.

\subsection{Skiving Stock Problem}
In the Skiving Stock Problem (SSP), we are given an integer $W$ and a set $I = \{1, \ldots, n\}$ of items, where each item $i \in I$ has a frequency $f_i \in \mathbb{Z}+$ and a size $w_i \in \mathbb{Z}+$, where $w_i < W$. The goal is to maximize the number of stock rolls with a size of at least $W$ while using at most $f_i$ of each item $I$. This problem is usually called the dual problem for the CSP and can be formulated using the SPF\@. The pricing is analogous to the CSP, where a pattern $P$ represents a subset of items satisfying ${W \leq \sum_{i = 1}^n a_i^{P} w_i \leq W_{\max} = 2 \cdot (W - 1)}$. Note that patterns with occupancy exceeding $W_{\max}$ do not need to be considered, as at least one item can be removed from them while maintaining an occupancy of at least $W$.

All properties presented for the CSP can be dualized to SSP, encompassing Waste Optimizations that enable us to enhance $W_{\max}$ to a value smaller than $2 \cdot (W - 1)$. For SSP, we implemented all features presented in the algorithm for CSP, except cutting planes. Moreover, as the CRF inequality has the same direction as the item inequalities, patterns belonging to this inequality may exhibit a negative reduced cost if its dual variable is disregarded. Hence, we prohibit our recovery algorithm from generating them.

Furthermore, we use the greedy heuristic proposed by~\citet{Peeters_2006} to produce the initial incumbent and initial columns for RLM and replace the BFD with it in the Rounding heuristic. This greedy heuristic consists of creating the patterns one by one and packing in the current pattern the smallest item $i$ that makes this pattern feasible or the greater remaining item if item $i$ does not exist.

For the computational experiments, we use the benchmarks A1, A2, A3, and~B generated by~\citet{Martinovic_2020}, and all benchmarks from BPP Lib, except Random class, which were extended by~\citet{Korbacher_2023} generating new instances for GI class with 750 and 1000 items of different sizes. We compare our algorithm with state-of-the-art algorithms MDISS~(\citet{Martinovic_2020}), NF-F~(\citet{Loti_2022}), and KIMS~(\citet{Korbacher_2023}), where all three algorithms use models based on AFF. The results for the MDISS algorithm were obtained by an AMD A10-5800K with 16 GB RAM (STPI $1491$), for the NF-F algorithm were obtained by an Intel Xeon E3-1245 v5 at 3.50GHz and 32 GB RAM (STPI $2249$), and for the KIMS were obtained by an i7-5930k CPU clocked at 3.5 GHz and 64 GB of RAM~(STPI $2050$)\@.

\begin{table}[t!]
\centering
\caption{Comparison with state-of-the-art algorithms for the SSP (time limit of 3600s).}
\label{Table6}
\footnotesize
\begin{tabular}{lrrrrrrrrrrrrrr}
\hline
               &                &                & \multicolumn{2}{c}{\textbf{MDISS}} &                & \multicolumn{2}{c}{\textbf{NF-F}} &                &  \multicolumn{2}{c}{\textbf{KIMS}}  &                & \multicolumn{3}{c}{\textbf{Our}}               \\
\cline{4-5} \cline{7-8} \cline{10-11} \cline{13-15} 
\textbf{Class} & \textbf{Total} & \hspace{0.15cm} & \textbf{~Opt~}   & \textbf{Time}   & \hspace{0.15cm} & \textbf{~Opt~}   & \textbf{Time}   & \hspace{0.15cm}& \textbf{~Opt~}   & \textbf{Time}  & \hspace{0.15cm} & \textbf{~Opt~} & \textbf{Time} & \textbf{Cols} \\
\midrule
A1 & 1260       && \textbf{1260} & 0.1  && \textbf{1260} & 0.1  && – & –                && \textbf{1260} & \textbf{0.005} & 38.8 \\
A2 & 1050       && 1011 & 250.9         && 1024 & 148.7         && 1023 & 137.3         && \textbf{1050} & \textbf{1.1} & 968.1 \\
A3 & 600        && – & –                && – & –                && 575 & 171.2          && \textbf{600} & \textbf{0.3} & 348.9 \\
AI202 & 50      && \textbf{50} & 26.3   && – & –                && \textbf{50} & 22.5   && \textbf{50} & \textbf{0.5} & 1104.9 \\
AI403 & 50      && 31 & 2016.9          && – & –                && 28 & 2024.4          && \textbf{50} & \textbf{8.9} & 3422.4 \\
AI601 & 50      && – & –                && – & –                && 1 & 3588.2           && \textbf{48} & \textbf{214.9} & 6924.4 \\
AI802 & 50      && – & –                && – & –                && – & –                && \textbf{44} & \textbf{552.8} & 11544.2 \\
AI1003 & 50     && – & –                && – & –                && – & –                && \textbf{40} & \textbf{866.2} & 16882.2 \\
ANI201 & 50     && \textbf{50} & 238.5  && – & –                && \textbf{50} & 274.5  && \textbf{50} & \textbf{0.7} & 1173.5 \\
ANI402 & 50     && 0 & 3600.0           && – & –                && 1 & 3586.0           && \textbf{50} & \textbf{5.7} & 3366.7 \\
ANI600 & 50     && – & –                && – & –                && – & –                && \textbf{49} & \textbf{99.3} & 6354.3 \\
ANI801 & 50     && – & –                && – & –                && – & –                && \textbf{46} & \textbf{398.0} & 11111.8 \\
ANI1002 & 50    && – & –                && – & –                && – & –                && \textbf{41} & \textbf{853.3} & 16108.5 \\
B & 160         && 126& 970.7           && \textbf{160} & 61.4  && \textbf{160} & 38.7  && \textbf{160} & \textbf{2.9} & 1596.4 \\
FalkenauerT& 80 && \textbf{80} & 0.9    && – & –                && – & –                && \textbf{80} & \textbf{0.5} & 1771.7 \\
FalkenauerU& 80 && \textbf{80} & 0.1    && – & –                && – & –                && \textbf{80} & \textbf{0.02} & 121.6 \\
GI125 & 80      && 40 & 1802.1          && – & –                && \textbf{80} & 19.4   && \textbf{80} & \textbf{14.8} & 419.8 \\
GI250 & 80      && 40 & 1854.0          && – & –                && \textbf{80} & 109.6  && \textbf{80} & \textbf{34.6} & 1120.8 \\
GI500 & 80      && 2 & 3584.2           && – & –                && 79 & 596.9           && \textbf{80} & \textbf{95.2} & 2866.0 \\
GI750 & 40      && – & –                && – & –                && 39 & 664.2           && \textbf{40} & \textbf{219.1} & 5528.9 \\
GI1000 & 40     && – & –                && – & –                && 38 & 1570.3          && \textbf{39} & \textbf{397.9} & 8461.9 \\
Hard & 28       && \textbf{28} & 3.0    && – & –                && – & –                && \textbf{28} & \textbf{0.1} & 435.4 \\
Scholl & 1210   && \textbf{1210} & 8.8  && – & –                && – & –                && \textbf{1210} & \textbf{0.2} & 221.6 \\
Schwerin & 200  && \textbf{200} & 0.2   && – & –                && – & –                && \textbf{200} & \textbf{0.04} & 152.3 \\
Waescher & 17   && \textbf{17} & 253.3  && – & –                && – & –                && \textbf{17} & \textbf{0.001} & 2.8 \\
\bottomrule
\end{tabular}
\end{table}

Table~\ref{Table6} showcases average run times and the number of instances optimally solved by each algorithm. While there is no widespread agreement on the instances used by various authors, we applied our framework across all of them, including AI and ANI instances with 600 items or more. Our algorithm surpasses all other state-of-the-art algorithms, marking the first occasion of solving all instances from classes A2 and A3. In SSP, the classes AI and ANI also involve checking if there exists a solution that comprises only patterns with occupancy equal to $W$, making the difficulty of solving these instances the same in both CSP and SSP\@. Although our CSP solver left 13 instances unsolved in these classes, our SSP solver left 32 instances unsolved. This discrepancy could result from the absence of cutting planes in the SSP solver or the pricing problem being more costly at the root node. Lastly, the initial heuristic and volume bound solved 16 out of 17 Waescher instances, indicating that instances in SSP may be notably easier than in CSP\@.

\subsection{Identical Parallel Machines Scheduling with Minimum Makespan}

In Identical Parallel Machines Scheduling~(IPMS), we need to assign a set $J$ of jobs to a set of $m$ identical machines. Preemption is not allowed, and each machine processes just one job at a time. Each job $j \in J$ has a processing time $s_j \in \mathbb{Z}_+$. A possible objective is minimizing the makespan, that is, the completion time of the last job to finish, in which case the problem is denoted by $P_m  \vert  \vert  C_{\max}$.
Also, note that $P_m  \vert  \vert  C_{\max}$ has a strong relation with BPP, since a solution with value at most $W$ for
$P_m  \vert  \vert  C_{\max}$ is a solution with at most $m$ bins for the BPP with bin capacity $W$ and vice-versa.

In the literature, $P_m  \vert  \vert  C_{\max}$ was explored by~\citet{DellAmico_2008},~\citet{Mrad_2018}, and~\citet{Gharbi_2022}. In~\citet{DellAmico_2008}, the authors use this relation with the BPP, using an exact algorithm for BPP based on SCF with CG and a binary search to solve $P_m  \vert  \vert  C_{\max}$. In~\citet{Mrad_2018}, the authors use a formulation based on the arc-flow formulation for the BPP, and~\citet{Gharbi_2022} improves this formulation using graph compression and uses the meta-heuristic Variable Neighborhood Search to find reasonable initial solutions.

In this work, we adopt the approach from the first paper, addressing $P_m  \vert  \vert  C_{\max}$ through a binary search combined with our CSP solver. Our binary search~(BS) algorithm comprises two stages, with its pseudocode provided in Appendix~\ref{app::IPMS}. The first stage invokes our CSP solver, progressively augmenting the current lower bound~(LB) with a sequence of powers of 2 until finding a solution that improves the current upper bound~(UB). The second stage tests the midpoint between LB and UB as the usual BS\@. This first stage is employed because we think that an optimal solution is more likely to be closer to LB than to UB\@. We initialize LB with the volume bounds $\max(\lceil \sum_{i \in I} \frac{d_i w_i}{m} \rceil, \max_{i \in I} w_i)$, and set UB based on the solution provided by the Longest Processing Time heuristic, a $\frac{4}{3}-$approximation algorithm \citep{Graham_1969}.

We also leverage the information acquired by the solver between optimizations, initializing the RLM in each iteration with all valid columns from the previous iteration. This strategy yields a beneficial warm start, particularly when the bin capacity exceeds that of the previous optimization. Lastly, in the CSP solver, we set an upper bound of $m + 1$, halting optimization upon discovering an integer solution of value less than or equal to $m$ (or if it is proven to not exist).

\begin{table}[b!]
\centering
\caption{Results for the \boldsymbol{$P_m || C_\text{max}$} using a time limit of 1200s.}
\label{Table7}
\footnotesize
\begin{tabular}{cccrrrrrrrr}
\toprule
                                   &                                       &                & \multicolumn{2}{c}{\textbf{M\&S}} & \multicolumn{2}{c}{\textbf{G\&B}} &\multicolumn{4}{c}{\textbf{Our}}                            \\ \cmidrule(lr){4-5}\cmidrule(lr){6-7}\cmidrule(lr){8-11}
\multicolumn{1}{c}{\textbf{n/m}} & \multicolumn{1}{c}{\textbf{Classes}} & \textbf{Total} & \textbf{Opt}    & \textbf{Time} & \textbf{Opt}    & \textbf{Time}    & \textbf{Opt} & \textbf{Time} & \textbf{Cols} & \textbf{Cuts} \\
    \midrule
       2 &     1-5 &    500 &   \textbf{500} &   0.4 &   \textbf{500} &  0.2 &  \textbf{500} & \textbf{0.005} &   28.1 &  0.0 \\
       2 &     6-7 &    200 &   \textbf{200} &   1.8 &   \textbf{200} &  0.6 &  \textbf{200} & \textbf{0.02} &  123.6 &  0.01 \\
    2.25 &     1-5 &    500 &   \textbf{500} &   3.7 &   \textbf{500} &  0.8 &  \textbf{500} & \textbf{0.03} &  218.6 &  0.2 \\
    2.25 &     6-7 &    200 &   \textbf{200} & 152.5 &   \textbf{200} & 10.0 &  \textbf{200} & \textbf{0.3} &  1016.5 &  7.3 \\
     2.5 &     1-5 &    500 &   \textbf{500} &   4.2 &   \textbf{500} &  0.8 &  \textbf{500} & \textbf{0.02} &  124.2 &  0.04 \\
     2.5 &     6-7 &    200 &   183 & 249.0 &   \textbf{200} & 25.7 &  \textbf{200} & \textbf{0.1} & 544.2 &  3.8 \\
    2.75 &     1-5 &    500 &   \textbf{500} &   5.6 &   \textbf{500} &  1.3 &  \textbf{500} & \textbf{0.04} & 284.7 &  0.4 \\
    2.75 &     6-7 &    200 &   154 & 415.7 &   199 & 60.0 &  \textbf{200} & \textbf{0.4} &  1250.2 &  7.8 \\
       3 &     1-5 &    500 &   \textbf{500} &   2.6 &   \textbf{500} &  1.5 &  \textbf{500} & \textbf{0.03} &  226.0 &  0.7 \\
       3 &     6-7 &    200 &   186 & 207.5 &   \textbf{200} & 28.9 &  \textbf{200} & \textbf{0.2} &  715.4 &  6.4 \\
    \bottomrule
\end{tabular}
\end{table}

Table~\ref{Table7} compares our algorithm with state-of-the-art algorithms proposed by~\citet{Mrad_2018} and~\citet{Gharbi_2022} for the $P_m  \vert  \vert  C_{\max}$, all of them using a time limit of $1200$ seconds. We obtain the results for both algorithms from~\citet{Gharbi_2022}, which run in an Intel Core i7-4930k (3.4 GHz) and 34 GB RAM (STPI $1971$), and CPLEX Solver 12.10. Due to space constraints, for each proportion $\frac{n}{m}$, we group classes 1 to 5 and classes 6 and 7 into the same cell.

The state-of-the-art algorithm by~\citet{Gharbi_2022} left only one unsolved instance, a significant improvement over~\citet{Mrad_2018}, who left 77 instances unsolved. However, our framework solves all instances in less than 1 second on average, making it dozens of times faster than~\citet{Gharbi_2022} across all classes. This performance aligns with expectations, given that these instances were randomly generated straightforwardly, a scenario where our framework excels in CSP\@.

\subsection{Ordered Open-End Bin Packing Problem}
In the Ordered Open-End Bin Packing Problem~(OOEBPP), we are given a positive integer $W$ and a set $I = \{1, \ldots, n\}$ of items, where each item $i \in I$ has a weight $w_i$ and an arrival time $p_i = i$. In the arrival order, the objective is to pack $I$ into bins with capacity $W$ while minimizing the number of bins used. Moreover, the last item in a bin (named overflow item) can be only partially packed. This means the bin capacity is considered met if the total weight packed minus the weight of the item with the largest $p_i$ is less than $W$.

The two leading exact algorithms for the OOEBPP were introduced by \citet{Ceselli_2008} and \citet{Loti_2022}. The former presents a formulation grounded in set covering, defining a set of patterns $\mathcal{P}_i$ for each item $i \in I$, with $i$ as the overflow item. The latter leverages the aforementioned set covering model to derive an arc-flow model for OOEBPP\@. 

In the OOEBPP, we utilize all features of our framework except waste optimization. We then detail the adaptation process for the formulation, pricing problem, and branching scheme. The remaining features are implemented in the same manner as for the CSP\@. For simplicity, we represent each item as a triple $(w_i, y_i, p_i)$, where $y_i$ denotes the weight of the item that must fit inside the bin if item $i$ is the last one packed. The original items in the instance have $y_i = 1$, while superitems resulting from merging two (super)items $i$ and $j$, where $p_i < p_j$, are represented as triples $(w_i + w_j, w_i + y_j, p_j)$.

If we omit explicitly segregating the set of patterns by overflow items in the formulation proposed by \citet{Ceselli_2008}, the resulting formulation is analogous to the SCF for the CSP\@. Following this approach, the only distinction is the explicit consideration of each item as an overflow item in the pricing problem. To achieve this, we employ dynamic programming for the Knapsack Problem using items sorted in ascending order of $p_i$. Thus, for each $1 \leq i \leq |I|$, we attempt to identify violated patterns where $i$ is the last item packed, invoking the recover pattern algorithm (Algorithm~\ref{alg::enumeration}) for the state $(i - 1, W - y_i)$ with a partial pattern $\overline{P}$ comprising item $i$. To prevent a slow pricer, due to the fixed ordering of items, we restrict the number of simultaneous active cutting planes to 100 in RLM\@.

We employ the Best-Fit Decreasing-Time (BFDT) heuristic, introduced by \citet{Ceselli_2008}, both in the Rounding heuristic and for generating the initial incumbent and columns for RLM\@. In the BFDT heuristic, items are processed in descending order of arrival time. Each item $i$ is packed into the fullest bin $b$ where it fits, or a new bin is created if $b$ does not exist.

One possible branching approach is to utilize the Ryan-Foster scheme, selecting the pair of items $(i, j)$ that maximizes the lexicographical order of $(-r_{ij}, p_j - p_i, w_i + w_j)$, with $\delta_{ij} \notin \mathbb{Z}$ and $p_j > p_i$, where $r_{ij}$ represents the priority of the pair of items $(i, j)$ as defined in Section~\ref{sec::scheme}. However, we observe that incrementing our scheme by utilizing the hierarchical strategy employed by \citet{Ceselli_2008}, which first branches based on whether an item is the last item packed or not, leads to greater efficiency.

Let $\mathcal{P}_i$ be the set of patterns where $i$ is the last item packed, and $\sigma_i$ denote the usage of item $i$ as the last item, i.e., $\sigma_i = \sum_{P \in \mathcal{P}_{i}} \lambda_P$. If there exists any item $i$ such that $0 < \sigma_i < 1$, we proceed to branch on the item with a value closer to $1$. The left branch considers only patterns where $i$ is the last item, while the right branch considers only patterns where $i$ is not the last item. The pricer efficiently manages this branching scheme by excluding item $i$ from dynamic programming on the left branch and not considering item $i$ as the last item during pattern recovery on the right branch. If $\sigma_i = 0$ or $\sigma_i \geq 1$ (where $\sigma_i > 1$ indicates overpacking of item $i$) for all items, we resort to branching using the aforementioned Ryan-Foster scheme.

For our computational experiments, we utilize the same benchmarks as~\citet{Loti_2022}. The first seven instance sets are sourced from two-dimensional BPP variants adapted to OOEBPP by~\citet{Ceselli_2008}, and are accessible via 2DPackLib~\citep{Iori_2022}. The remaining sets consist of randomly generated instances proposed by~\citet{Loti_2022}, each containing up to 1000 items, and were obtained by us upon request to them. Originally, they only tested instances with 50, 100, and 200 items. As these smaller instances were easy for our framework, we conducted experiments on the entire instance set.

Table~\ref{Table8} presents a comparison between our algorithm and leading algorithms for the OOEBPP, all executed under a time limit of 3600 seconds. The results for CR from~\cite{Ceselli_2008} were obtained using a Pentium IV 1.6 GHz with 512 MB of RAM (STPI $225$), while those for Arc-flow and NF-F from~\cite{Loti_2022} were obtained using an Intel Xeon E3-1245 v5 at 3.50GHz with 32 GB RAM (STPI $2249$). Excluding the classes where~\cite{Loti_2022} report execution times rounded down to ``$<0.1$'', our algorithm outperforms state-of-the-art algorithms in all instances except the BENG class. Instances in this class consist of items with numerous items per bin, and the set covering formulation requires significant time to solve the root node. Thus, a potential approach to enhance runtime is to employ a better primal heuristic while solving the root node, but such optimization lies beyond the scope of this article.

\begin{table}[t!]
\centering
\caption{Results for the OOEBPP using a time limit of 3600s.}
\label{Table8}
\footnotesize
\begin{tabular}{lrrrrrrrrrrr}
    \toprule
          &       & \multicolumn{2}{c}{\textbf{CR}} & \multicolumn{2}{c}{\textbf{Arc-flow}} & \multicolumn{2}{c}{\textbf{NF-F}} & \multicolumn{4}{c}{\textbf{Our}}                           \\ \cmidrule(lr){3-4} \cmidrule(lr){5-6} \cmidrule(lr){7-8} \cmidrule(lr){9-12} 
\textbf{Class}     & \textbf{Total} & \textbf{Opt}   & \textbf{Time}  & \textbf{Opt}     & \textbf{Time}     & \textbf{Opt}    & \textbf{Time}   & \textbf{Opt} & \textbf{Time} & \textbf{Cols} & \textbf{Cuts}\\ \midrule 
BENG      & 10    & \textbf{10}    & \textbf{0.1}            & \textbf{10}      & 0.3               & \textbf{10}     & \textbf{0.1}             & \textbf{10}  & 1.3          & 1155.7     & 10.9           \\
CGCUT     & 3     & \textbf{3}     & 0.1            & \textbf{3}       & 0.3               & \textbf{3}      & 0.4             & \textbf{3}   & \textbf{0.02}          & 130.7              & 0.0           \\
CLASS     & 500   & 499   & 8.7            & \textbf{500}     & 11.3              & \textbf{500}    & 3.7             & \textbf{500} & \textbf{0.09}          & 295.6                       & 2.7         \\
GCUT 1-4  & 4     & \textbf{4}     & 0.6            & \textbf{4}       & 0.1               & \textbf{4}      & \textless 0.1   & \textbf{4}   & 0.008             & 79.8                     & 0.0             \\
GCUT 5-13 & 9     &–    &–             & \textbf{9}       & 1.3               & \textbf{9}      & 0.1             & \textbf{9}   & \textbf{0.006}             & 73.1                        & 0.0             \\
HT        & 9     & \textbf{9}     & 0.1            & \textbf{9}       & \textless 0.1     & \textbf{9}      & \textless 0.1   & \textbf{9}   & 0.01          & 53.1                        & 19.7           \\
NGCUT     & 12    & \textbf{12}    & 0.1            & \textbf{12}      & \textless 0.1     & \textbf{12}     & \textless 0.1   & \textbf{12}  & 0.004          & 33.3                       & 0.0           \\
Random50       & 480   &–             &–             & \textbf{480}     & 1.7               & \textbf{480}    & 1.3             & \textbf{480} & \textbf{0.03}          & 226.4             & 1.4          \\
Random100      & 480   &–             &–             & 479              & 73.4              & 479             & 38.7            & \textbf{480} & \textbf{0.2}           & 684.7             & 6.1          \\
Random200      & 480   &–             &–             & 426              & 863.3             & 466             & 168.9           & \textbf{480} & \textbf{2.5}           & 1770.6            & 11.6         \\
Random300      & 480   &–             &–             &–               &–                &–              &–              & \textbf{474} & \textbf{26.7}          & 3198.0                    & 47.2         \\
Random400      & 480   &–             &–             &–               &–                &–              &–              & \textbf{473} & \textbf{59.5}          & 4643.6                    & 78.5         \\
Random500      & 480   &–             &–             &–               &–                &–              &–              & \textbf{472} & \textbf{79.2}          & 9327.4                    & 98.0         \\
Random750      & 480   &–             &–             &–               &–                &–              &–              & \textbf{473} & \textbf{91.6}         & 9879.4                    & 46.8         \\
Random1000     & 480   &–             &–             &–               &–                &–              &–              & \textbf{464} & \textbf{228.3}         & 14631.5                   & 32.1         \\
\bottomrule        
\end{tabular}
\end{table}

Moreover, observe that our algorithm can solve the smaller Random instances dozens of times faster than algorithms proposed by~\citet{Loti_2022}. Our framework continues solving many instances for the non-tested instances by~\citet{Loti_2022}, though it left 36 unsolved instances.

\subsection{Class-Constrained Bin Packing Problem}
In the Class-Constrained Bin Packing Problem~(CCBPP), we have positive integers $C$, $Q$, and $W$, and a set $I = \{1, \ldots, n\}$ of items, where each item $i \in I$ has a weight $w_i$ and a class identifier $q_i \in \{1, \ldots, Q\}$. The objective is to pack $I$ in bins of capacity $W$ while minimizing the number of bins used and satisfying the constraint that each bin can have items of at most $C$ different classes.

The sole exact algorithm in the literature for the CCBPP is a branch-and-price algorithm using the SCF proposed by~\citet{Borges_2020}. We use all features of our framework for the CCBPP except Waste Optimizations and the CRF heuristic. The former is rarely beneficial, and the latter is unnecessary since the RF heuristic is sufficient to solve the literature's benchmarks. Next, we explain how to adapt the formulation, the pricing problem, and the branching scheme. The remainder features are used as in CSP\@.

The pricing problem for SCF, disregarding conflicts and cutting planes,  entails the Class-Constrained Knapsack Problem (CCKP), for which the authors provide two dynamic programming approaches. The first has a time complexity of $\Theta(nW + W^2CQ)$, while the second has a time complexity of $\Theta(nWC)$. Opting for efficiency, we utilize the second approach. However, we observed redundancies in the recurrence presented by the authors. Hence, we propose a revised recurrence as follows:
\begin{equation}
\resizebox{0.8\textwidth}{!}{%
  $f (i, r, c, u)=\begin{cases}
    1, & \text{if $i = 0$, $w = 0$, or $c = 0$}.\\
    f (i - 1, r, c, u \land e) & \text{if $i \geq 1$ and $r < w_{I[i]}$}.\\
    \min (f (i - 1, r, c, u \land e), \\ \hspace{0.75cm} f(i - 1, r - w_{I[i]}, c - \overline{u}, e) - \overline{\pi}_{I[i]}) & \text{if $i \geq 1$ and $r \geq w_{I[i]}$},
  \end{cases}\label{ineq::dp_CCBPP}$
}
\end{equation}
where items in $I$ are sorted in non-decreasing order of $q_i$ and $i$ is the index of the current item. The variables $r$ and $c$ represent, respectively, the remaining capacity and number of classes. Additionally, $u$ and its complement $\overline{u}$ are binary values, where $u$ indicates if $q_{I[i]}$ was already discounted from $c$, while $e$ indicates whether $q_{I[i]}$ and $q_{I[i - 1]}$ are equal. The lower bound for $c_P(\vpi, \vrho, 1)$ corresponds to the value of $f(n, W, C, 0)$.

We adopt the Ryan-Foster scheme for this problem, selecting the pair of items $(i, j)$ to branch based on the lexicographical order of $(-r_{ij}, e_{ij}, w_i + w_j)$, with $\delta_{ij} \notin \mathbb{Z}$, where $r_{ij}$ represents the priority of the pair $(i, j)$ as defined in Section~\ref{sec::scheme}, and $e_{ij}$ equals 1 if $q_i = q_j$, and 0 otherwise. This approach prevents the creation of items with multiple classes and demonstrates strong performance across the literature's benchmarks.

Adapting Algorithm~\ref{alg::enumeration} for the CCBPP is straightforward. To accommodate items with multiple classes (which arise when merging two items), we ensure that the items are sorted in non-increasing order of class identifier. Then, we utilize the recurrence from CCKP, computed solely with the highest class identifier of each item, as a lower bound. Consequently, a class $q$ that is not the highest identifier of a multiple class item is discounted from $c$ only upon encountering the first item with a class identifier less than or equal to $q$. Additionally, we verify the fulfillment of the class constraint in base cases with $i > 0$.

Furthermore, we use the BFD heuristic in the Rounding heuristic as well as for producing the initial incumbent and initial columns for RLM\@. In other words, we process the items in non-increasing order of weight and pack the current item $i$ in the fullest bin $b$ that $i$ can be packed satisfying its capacity and class constraint or create a new one if $b$ does not exist.

\begin{table}[t!]
\caption{Comparison between our algorithm and the state-of-the-art for the CCBPP.}
\label{Table9}
\footnotesize
\centering
\begin{tabular}{ccccrrrrrrr}
\toprule
\multicolumn{2}{c}{\textbf{Class}} &                                                              & & \multicolumn{3}{c}{\textbf{\citet{Borges_2020}}}                                                & \multicolumn{4}{c}{\textbf{Our}}                                                                              \\ \cmidrule(lr){1-2} \cmidrule(lr){5-7} \cmidrule(lr){8-11}
{\textbf{Q}} & {\textbf{C}} & \textbf{Total} & \textbf{Non-IRUP \footnotemark[1]} & {\textbf{Opt}}          & {\textbf{Time}}   & \textbf{Cols}    & {\textbf{Opt}}          & {\textbf{Time}}  & {\textbf{Cols}}    & \textbf{Cuts}  \\ \midrule
10 & 2 & 120 & 0 & \textbf{120} & 4.3 & 833.3 & \textbf{120} & \textbf{0.6} & 1019.7 & 1.2 \\
10 & 3 & 120 & 0 & \textbf{120} & 3.1 & 739.9 & \textbf{120} & \textbf{0.4} & 954.9 & 0.0 \\
10 & $\chi$ & 120 & 0 & 97 & 178.1 & 722.7 & \textbf{120} & \textbf{0.1} & 565.8 & 0.0 \\
25 & 2 & 120 & 0 & \textbf{120} & 10.0 & 1040.6 & \textbf{120} & \textbf{1.1} & 1014.1 & 5.6 \\
25 & 3 & 120 & 0 & \textbf{120} & 2.4 & 702.3 & \textbf{120} & \textbf{0.4} & 1017.2 & 0.0 \\
25 & $\chi$ & 120 & 0 & 110 & 87.6 & 697.5 & \textbf{120} & \textbf{0.2} & 695.5 & 0.0 \\
50 & 2 & 120 & 22 & 96 & 185.4 & 2871.4 & \textbf{120} & \textbf{6.6} & 1420.8 & 33.4 \\
50 & 3 & 120 & 0 & \textbf{120} & 2.1 & 725.2 & \textbf{120} & \textbf{0.4} & 1048.1 & 0.0 \\
50 & $\chi$ & 120 & 1 & 119 & 13.3 & 790.1 & \textbf{120} & \textbf{0.2} & 681.7 & 0.7 \\
\bottomrule
\end{tabular}
\footnotetext[1]{Amount of non-IRUP instances in the class} 
\end{table}

For the computational experiments, we adopt the benchmark proposed by~\citet{Borges_2020}, which feature instances with $n = 200$ items, bin capacities $W \in \{100, 150, 200\}$, number of classes $Q \in \{10, 25, 50\}$, and number of the class by bin $C \in \{2, 3, \chi \}$, where $\chi = \lfloor n \mathbin{/} {\lfloor \sum_{i \in I} w_i} \mathbin{/} W \rfloor \rfloor$. This benchmark contains instances with $C = 1$, which are unions of BPP instances, one for each class $q \in {1, \ldots, Q}$, as noted by the authors. While solving them individually using a BPP solver would be more efficient, we follow~\citet{Borges_2020} and treat them as CCBPP instances. Notably, in instances with $C = 1$, the lower bound $\lceil \sum_{P \in \mathcal{P}}\vlambda_P \rceil$ can be weak, since each relaxation of BPP instance can have a gap for the next integer almost equal to 1, leading to a total gap almost equal to $Q$. Seen this, when $C = 1$, our safe bound is computed using the aggregated lower bound $\sum_{q = 1}^Q \lceil \sum_{P \in \mathcal{P}_q}\vlambda_P \rceil$, where $\mathcal{P}_q$ is the set of patterns $P$ with items of the class $q$.

The data provided in~\citet{Borges_2020} does not allow an easy comparison with our algorithm. To facilitate a fair evaluation, we re-ran their algorithm using the same resources as ours (Gurobi and GCC versions, and hardware). Table~\ref{Table9} presents our algorithm with the best version of the algorithm proposed by~\citet{Borges_2020}. Notably, our algorithm outperforms the current state-of-the-art, solving all instances at least five times faster. Additionally, we noticed that among the unsolved instances left by~\citet{Borges_2020}, the classical lower bound of rounding up the relaxation value proves very weak in roughly half of them. Hence, leveraging our aggregated lower bound might enable their algorithm to tackle these instances.

Furthermore, 22 instances are non-IRUP among the 120 instances with $Q = 50$ and $C = 2$. These cases exhibit an average discrepancy of $1.38$ between the solution value obtained by RLM without cutting planes and the optimal solution value. Notably, one instance demonstrates a gap exceeding $2$, which contrasts with the other problems studied in this article by satisfying the MIRUP conjecture\@. However, upon incorporating SRIs, these gaps diminish to nearly $1$. This improvement comes with the trade-off of introducing hundreds of additional cuts, substantially inflating the time required to solve each relaxation node. This increased computational overhead is exacerbated by the need to explore numerous nodes to close the remaining gap.

\section{Conclusions\label{sec::Conclusions}}
In this work, we propose a framework leveraging the SCF/SPF to tackle five cutting stock problems already studied in the literature. Our key contribution lies in a CRG process employing a multi-pattern generation and diversification strategy, ensuring convergence in a few iterations and outperforming existing algorithms significantly. Notably, our advancements are most pronounced in SSP, $P_m \vert \vert C_{\max}$, CCBPP, and OOEBPP, where our framework achieves success in solving instances previously unaddressed in the literature. Additionally, our primal heuristics prove highly effective, often solving many instances in the root node or exploring only a few nodes. Furthermore, we exploit the CSP symmetry and employ strategies for produce a leaner model, enabling successful solution of the challenging AI and ANI instance classes.

In future research, our insights could yield promising outcomes for other problems featuring the SCF/SPF, strong relaxations, and pseudo-polynomial pricers. Moreover, the success of our algorithm with AI and ANI instances stems from the proximity of the relaxation polyhedra to the convex hull of integer solutions. Thus, we introduced a new set of instances lacking this property, which emerged as a challenging benchmark for the CSP\@. This novel benchmark's structure can be extended across various problem domains beyond cutting and stock problems, potentially leading to a new research line.

\section*{Acknowledgments}
Supported by the São Paulo Research Foundation (FAPESP) grants \#2015/11937-9 and \#2022/05803-3; the Brazilian National Council For Scientific and Technological Development (CNPq) grant \#312345/2023-2; and Teaching, Research and Extension Support Fund of University of Campinas (FAEPEX/UNICAMP) grant \#2372/23.

\begin{appendices}

\section{Proof of the Conflict Propagation's Lemma\label{app::propagation}}
Let $v(i, j)$ denote the current B\&B node, where the Extended Ryan-Foster Scheme selects the pair  $(i, j)$ of superitems for branching. Let $k$ represent the item such that $w_k = w_i + w_j$. Additionally, node $v(i, j)$ possesses a set of superitems $\overline{I}$, where each superitem $i$ is a nonempty subset $S \subseteq I$ (the items from the original instance). Alongside, there exists a conflict graph $G = (\overline{I}, E)$, where each edge $\{k, l\} \in E$ represents a conflict introduced by a right branch of our branching scheme or by conflict propagation. We denote by $E_i$ the set of superitems that have a conflict with superitem $i$ in~$G$, i.e., the adjacency list of vertex $i$ in $G$.

An alternative approach to Definition~\ref{definition::scheme} is to consider the edges added to $G$ rather than focusing solely on modifications in the conflict lists. Let $G = (\overline{I}, E)$ and $G' = (\overline{I}', E')$ represent the conflict graphs before and after branching at the current node $v(i, j)$. In this context, the left branch defines 
$E' = E \cup \{\{a, k\} \colon \{a, i\} \in E \vee \{a, j\} \in E\}$,
while the right branch defines 
$E' = E \cup \{\{i, j\}\}$. An interesting observation is that the conflict added to a node $v(i, j)$ remains valid for all its descendants, regardless of whether the item demand becomes zero at some point.

We define an \emph{improvement solution} as one better than the incumbent solution. Subsequently, we prove the feasibility of maintaining only one conflict list for each item size, as stated in Lemma~\ref{lemma::propagation}, addressing conflict propagation as presented in Section~\ref{sec::conflicts}. 

If there is an item $k' \in \overline{I}$ with $w_{k'} = w_i + w_j$, a priori, we cannot group $k$ and $k'$ into one item by summing their demands, as it could lead to a loss of an improvement solution. The same is valid for an item $i'$ with $w_{i'} = w_i$ or $j'$ with $w_{j'} = w_j$. 

Next, we prove that this grouping is valid through some claims.

\begin{claim}
    It is valid to define $E_k' = E_{k'} \cup E_i \cup E_j$ on the left branch of node $v(i, j)$.
\end{claim}

\emph{Proof: } It is straightforward that $E_i \cup E_j \subseteq E_k'$, as $k$ is created by merging $i$ and $j$, thus we show it is valid to define $E_{k'} \subseteq E_k'$. This is untrue only if there exists $h \in E_{k'}$ such that there is an improvement solution $\vlambda$ in the subtree rooted at $v(i, j)$ containing a pattern $P_1$ with $h$ and $k$ together. 

Let $v(k', h)$ be the ancestor of $v(i, j)$ that added the conflict $\{k', h\}$ in its right branch. Let $P_2$ be a pattern containing item $k'$ in $\vlambda$, and let $\vlambda'$ be the equivalent solution obtained by swapping $k'$ and $k$, with $P_1' = P_1 \cup \{k'\} \setminus \{k\}$ and $P_2' = P_2 \cup \{k\} \setminus \{k'\}$. It is possible that \( k' \) does not exist as an individual item because it was merged with other items to form a superitem \( l \). In this case, we redefine \( l \) to be created using \( k \) instead of \( k' \). If there is no conflict in $P_1'$ and $P_2'$, then $\vlambda'$ is an improvement solution in the search space of the left child of $v(k', h)$, which leads to a contradiction (as we prioritize left in the DFS order).

Suppose, on the contrary, that either there is a conflict $\{k', g\}$ in $P_1'$, with $g \in E_{k'}$, or there is a conflict $\{i, g\}$ or $\{j, g\}$ in $P_2'$, with $g \in E_i \cup E_j$. The conflicts $\{i, g\}$ and $\{j, g\}$ might occur inclusive within the superitem $l$, but it does not affect our argument. Without loss of generality, as all cases are analogous, let $\{k', g\}$ be the earliest conflict. Observe that this conflict was introduced by the ancestor $v(k', g)$ of $v(k', h)$ in its right branch. Thus,  $\vlambda'$ is a feasible and improvement solution in the left child of $v(k', g)$, which is, again, a contradiction.
~\hfill$\blacksquare$

\begin{claim}
It is valid to define $E_{k'}' = E_{k'} \cup E_i \cup E_j$ on the left branch of the current node $v(i, j)$.
\end{claim}

\emph{Proof:} Clearly, $ E_{k'} \subseteq E_{k'}'$, as by the branching scheme, conflicts are inherited. Now, it is invalid to define $E_i \cup E_j \subseteq E_{k'}'$ only if there exists an item $h \in E_i \cup E_j$ such that there is an improvement solution $\vlambda$ in the subtree of $v(i, j)$ containing a pattern $P_1$ with $h$ and $k'$ together. As said before, if \( k' \) was used to create a superitem $l$ at $v(i, j)$, redefine $l$ to be created with \( k \).

Without loss of generality, assume that $h \in E_i$ and that $v(i, h)$ is the earliest ancestor of $v(i, j)$ with $h \in P_1 \cap E_i$. Let $\vlambda'$ be an equivalent solution obtained by swapping $k'$ and $k$. Analogously to the previous claim, we can define $P_2$, $P_1'$, and $P_2'$ and argue that $\vlambda'$ is an improvement solution in the left child of $v(i, h)$ if it is feasible, or in the left child of an ancestor of $v(i, h)$ otherwise, which leads to a contradiction. Thus, it is valid to define $E_{k'}' = E_{k'} \cup E_i \cup E_j$.~\hfill$\blacksquare$

\begin{claim}
    It is valid to define $E_{i'}' = E_{i}' = E_{i} \cup \{j\}$ and $E_{j'}' = E_{j}' = E_{j} \cup \{i\}$ on the right branch of the current node $v(i, j)$.
\end{claim}

\emph{Proof:} 
We introduce a conflict between $i$ and $j$ on the right branch, updating $E_i' = E_i \cup \{j\}$ and $E_j' = E_j \cup \{i\}$. Suppose that we cannot define $E_{i'}'$ equal to $E_i'$ (the same argument applies to $j$ and $j'$). Using a similar argument than before, this is because there is an improvement solution where item $i'$ is with an item $g \in E_i$. By choosing a $g$ such that $v(i, g)$ is the earliest ancestor of $v(i, j)$, we can swap $i$ and $i'$ to obtain an improvement solution in the left child of $v(i, g)$, which leads to a contradiction.
~\hfill$\blacksquare$

Since the conflict lists of all items with the same weight are identical, we can combine them into one item by summing their demands, and this implies the correctness of Lemma~\ref{lemma::propagation}.

\section{Waste Optimization at Root\label{app::waste}}

During the column generation at the root node, one iteration can generate patterns with a large waste of space. Over subsequent iterations, these patterns cannot belong to any optimal RLM solution. As these patterns become obsolete, removing them from the model can accelerate the algorithm. Therefore, we propose a technique to shrink the RLM at the root.

Let~$z(\vlambda)$ be the value of the solution $\vlambda$ and $\Ws = \sum_{i \in I} d_i w_i$ be the total size of the items. We use the following rule to shrink the RLM at the root node.

\begin{Rule}\label{rule::root}
Given an intermediary solution $\Zl$ in the root-solving process, add the constraint that all patterns in the primal model or future generated have a waste of space at most $\Rr = z(\Zl) \cdot W - \Ws$.
\end{Rule}

Note that this rule also prevents the generation of new patterns with waste greater than $\Rr$ in the future, thereby avoiding the addition of obsolete patterns to the model. Next, we demonstrate that by applying Rule~\ref{rule::root}, we can obtain a valid lower bound for the problem under some conditions.

\begin{lemma}\label{lemma:root}
Let $z^\text{safe}_R$ be the safe lower bound obtained after we apply Rule~\ref{rule::root} with an intermediate solution $\Zl$. If $\lceil z^\text{safe}_R \rceil  \leq \lceil z(\Zl) \rceil$, then $z^\text{safe}_R$ is a valid lower bound for the problem, i.e., $\lceil {z^\text{safe}_R \rceil \leq z(\Zi)}$, where $\Zi$ is an optimal integer solution.
\end{lemma}
\textbf{Proof:} If all used patterns of $\Zi$ have waste less than or equal to $\Rr$, then $\lceil z^\text{safe}_R \rceil \leq z(\Zi)$, since $\Zr^*$ can use the same patterns of LM\@. Thus, suppose that there is a pattern $P$ with waste $R' > \Rr$ in  $\Zi$ with ${(\Zi)_P \geq 1}$. In this case, note that ${z(\Zi) \cdot W > \Ws + \Rr}$ and, thus, ${\Rr < z(\Zi) \cdot W - \Ws}$. By definition, it follows that ${z(\Zl) \cdot W - \Ws < z(\Zi) \cdot W - \Ws}$, and ${z(\Zl) < z(\Zi)}$. 
As $z(\Zi)$ is an integer, we conclude that ${\lceil z^\text{safe}_R \rceil \leq \lceil z(\Zl) \rceil \leq z(\Zi)}$.~\hfill$\blacksquare$

It is worth noting that Lemma~\ref{lemma:root} does not address the case where $\lceil z^\text{safe}_R \rceil > \lceil z(\Zl) \rceil$. However, in our experiments, we did not observe any instances where this inequality was violated. If it does occur, we can ensure the algorithm's correctness by removing the waste constraint and continuing the column generation process. Additionally, to ensure that our pricer generates only patterns with waste up to $\Rr$, we can modify the base case $(i = 0)$ of dynamic programming presented in Section~\ref{sec::recover}, setting $dp(0, r) = K$ if $r \leq \Rr$, and $dp(0, r) = \infty$ otherwise.

\section{Feature tests in CSP\label{app::features_tests}}

In this section, we present tests involving the techniques discussed in this paper. We evaluate our final algorithm, which uses all features, and eleven different versions. Each of the first nine versions excludes one of the following features from the final algorithm: multiple pattern generation; the RF heuristic; RF and CRF heuristics; the splay operation; the historic branching; splay operation and historic branching; the smaller tolerance $\epsilon' = 2.5 \cdot 10^{-12}$ introduced in Section\ref{sec::num_safe}; the dual inequalities and binary pricing; and the Model-Cleaning-by-Reduced-Cost (MCRC). Additionally, there is one version without conflict propagation and grouping of items, thus functioning as a BPP solver instead of a CSP solver (where an item $i$ with demand $d_i$ is considered to be $d_i$ independent items). The last version employs tan extension of \emph{standard branching scheme} used by \citet{Wei_2020} and \citet{Baldacci_2024}, which does not use the splay operation or historic branching, and branches on the pair of items $(i, j)$ where $\delta_{ij}$ is the farthest to an integer value. This is an extension because both works considers the BPP and, thus, they choose the pair of items $(i, j)$ with~$\delta_{ij}$ closer to $0.5$.

Regarding the modifications in the algorithm, we use our pricing structure in the version without multiple pattern generation restricted to return only the smallest reduce-cost pattern. The version that removes RF must also remove CRF since the second one depends on the first one. In the version that turns off the historic branching, the priority $r_{ij}$ is set to zero for all pairs of items. Due to a small redundancy between the splay operation and historic branching, we also have a version that removes both. We use the standard tolerance $\epsilon = 10^{-9}$ in the version without the smaller tolerance.

Tables~\ref{Table2},~\ref{Table3}, and~\ref{Table3b} present the results for this comparison using a time limit of 600 seconds, showing the difference in the number of instances solved (in absolute value) and the average time (in percentage) between each one of the eleven versions and the complete version.

The \textsc{No Multiple Patterns} version shows that multiple pattern generation is crucial. Without it, no AI and ANI instances with more than 600 items were solved, and the number of solved GI instances was reduced by 96. The \textsc{No CRF} and \textsc{No RF No CRF} versions highlight the importance of RF and especially CRF, increasing the number of solved AI instances by 23.

Disabling the splay operation increased the time required to solve the Random and Hard classes by 60\% and 35\%, respectively. Historic branching is essential for solving all Hard and Random instances and is necessary for solving ANI instances with 600 or more items. When both features are removed, the number of unsolved Hard, Random, and ANI instances increases even more, indicating an overlap in their effects. In the \textsc{Standard Branching} version, the number of unsolved instances in the Random and Hard classes decreases compared to the \textsc{No Splay No Historic} version. However, the number of unsolved instances and the time spent on the AI and ANI classes increase, even for instances with 400 items.

The \textsc{Standard EPS} version indicates that the smaller tolerance is essential for solving ANI instances with more than 600 items. The \textsc{No Dual Ineq. No B. Pricing} version shows that dual inequalities and binary pricing improve time performance in several instance classes, such as GI, Schwerin, and Waescher. The \textsc{No MCRC} version shows that our variable elimination by reduced cost is necessary for solving some AI and ANI instances and to improve the time performance on other classes.

Moreover, the \textsc{No Grouping} version has significantly worse performance, taking ten times longer in several classes, and it leads to a \emph{memory limit exceeded} error in the GI class with a bin capacity $W = 1.5 \cdot 10^6$, even on a machine with 64 GB of RAM\@. We tested the inclusion of conflict propagation in this version to check if it leads to better performance, but the results indicated no significant difference. Thus, the main benefit of the Extended Ryan-Foster Scheme is not the conflict propagation itself, but the fact that conflict propagation allows for the grouping of items. The grouping of items breaks many symmetries in the linear relaxation polyhedron and the pricing problem, as well as enabling the use of the powerful Ryan-Foster scheme in non-binary problems.

\begin{table}[th!]
\centering
\small
\caption{Comparison between the first three versions and the final version using a time limit of 600s.}
\label{Table2}
\begin{tabular}{lrrrrrrrr}
\toprule
& \multicolumn{2}{c}{\textsc{\mline{Final}{Version}}} &      \multicolumn{2}{c}{\textsc{\mline{No Multiple}{Patterns}}}  & \multicolumn{2}{c}{\textsc{No CRF}}   & \multicolumn{2}{c}{\textsc{\mline{No RF}{No CRF}}}\\
Name        &   Opt & Time    &$\Delta$ Opt & $\Delta$ Time&$\Delta$ Opt & $\Delta$ Time&$\Delta$ Opt & $\Delta$ Time\\
\midrule
AI202 & 50 & 0.41 & 0 & 697.6\% & 0 & 0.0\% & 0 & 97.6\% \\
AI403 & 50 & 6.01 & 0 & 1066.4\% & -2 & 443.6\% & -4 & 822.1\% \\
AI601 & 48 & 46.18 & 0 & 657.7\% & -7 & 154.2\% & -9 & 217.4\% \\
AI802 & 47 & 93.67 & -47 & 540.6\% & -7 & 46.9\% & -9 & 83.9\% \\
AI1003 & 38 & 190.60 & -38 & 214.9\% & -1 & 1.9\% & -1 & 11.8\% \\
ANI201 & 50 & 0.48 & 0 & 654.2\% & 0 & 2.1\% & 0 & -14.6\% \\
ANI402 & 50 & 2.69 & 0 & 2364.7\% & 0 & -3.3\% & 0 & -10.0\% \\
ANI600 & 50 & 13.42 & 0 & 2100.5\% & 0 & -2.1\% & 0 & -7.1\% \\
ANI801 & 49 & 38.78 & -49 & 1447.3\% & 0 & -4.6\% & 0 & -4.5\% \\
ANI1002 & 44 & 139.95 & -44 & 328.8\% & 0 & -3.1\% & 1 & -5.2\% \\
FalkenauerT & 80 & 0.10 & 0 & 230.0\% & 0 & 0.0\% & 0 & 360.0\% \\
FalkenauerU & 80 & 0.02 & 0 & 50.0\% & 0 & 0.0\% & 0 & 0.0\% \\
Hard & 28 & 7.05 & 0 & 8.5\% & 0 & -18.9\% & 0 & -52.1\% \\
GI AA & 60 & 18.25 & -20 & 1383.8\% & 0 & -3.2\% & -10 & 620.8\% \\
GI AB & 60 & 9.92 & -15 & 2806.9\% & 0 & -12.0\% & -5 & 567.0\% \\
GI BA & 59 & 71.88 & -38 & 521.6\% & 0 & -21.5\% & -10 & 136.8\% \\
GI BB & 60 & 33.07 & -23 & 1095.0\% & 0 & -9.3\% & -2 & 64.1\% \\
Random & 3840 & 0.05 & 0 & 240.0\% & 0 & 0.0\% & 0 & 280.0\% \\
Scholl & 1210 & 0.06 & 0 & 416.7\% & 0 & -16.7\% & 0 & 533.3\% \\
Schwerin & 200 & 0.02 & 0 & 50.0\% & 0 & 0.0\% & 0 & 50.0\% \\
Waescher & 17 & 0.31 & 0 & 54.8\% & 0 & 3.2\% & 0 & 183.9\% \\
\bottomrule
\end{tabular}
\end{table}
       
\begin{table}[H]
\centering
\small
\caption{Comparison between the splay / historic versions and final version using time limit of 600s.}
\label{Table3}
\centering
\begin{tabular}{lrrrrrrrrrr}
\toprule
& \multicolumn{2}{c}{\textsc{\mline{Final}{Version}}} &\multicolumn{2}{c}{\textsc{No Splay}} & \multicolumn{2}{c}{\textsc{No Historic}} & \multicolumn{2}{c}{\textsc{\mline{No Splay}{No Historic}}}  &\multicolumn{2}{c}{\textsc{\mline{Standard}{ Branching}}}\\
Name        &   Opt & Time    &$\Delta$ Opt & $\Delta$ Time&$\Delta$ Opt & $\Delta$ Time&$\Delta$ Opt & $\Delta$ Time&$\Delta$ Opt & $\Delta$ Time\\
\midrule
AI202 & 50 & 0.41 & 0 & 4.9\% & 0 & -2.4\% & 0 & 4.9\% & 0 & 70.7\% \\
AI403 & 50 & 6.01 & 0 & -37.4\% & 0 & 12.1\% & 0 & -10.6\% & -1 & 330.9\% \\
AI601 & 48 & 46.18 & 0 & 9.7\% & 0 & 13.8\% & 0 & 8.8\% & -3 & 138.1\% \\
AI802 & 47 & 93.67 & -1 & 14.7\% & -1 & -1.0\% & -1 & 3.9\% & -1 & 19.7\% \\
AI1003 & 38 & 190.60 & -1 & 2.9\% & -2 & 7.2\% & -1 & 7.0\% & -3 & 13.9\% \\
ANI201 & 50 & 0.48 & 0 & 0.0\% & 0 & -4.2\% & 0 & 0.0\% & 0 & 37.5\% \\
ANI402 & 50 & 2.69 & 0 & -3.7\% & 0 & 4.5\% & 0 & 3.7\% & -1 & 433.1\% \\
ANI600 & 50 & 13.42 & 0 & -1.3\% & 0 & 97.2\% & 0 & 96.3\% & -2 & 145.2\% \\
ANI801 & 49 & 38.78 & 0 & -2.5\% & -1 & 14.8\% & -1 & 15.8\% & -1 & 30.1\% \\
ANI1002 & 44 & 139.95 & 0 & -4.9\% & -2 & 12.7\% & -3 & 19.8\% & -3 & 10.6\% \\
FalkenauerT & 80 & 0.10 & 0 & 10.0\% & 0 & 0.0\% & 0 & 0.0\% & 0 & 0.0\% \\
FalkenauerU & 80 & 0.02 & 0 & 0.0\% & 0 & 0.0\% & 0 & 0.0\% & 0 & 0.0\% \\
Hard & 28 & 7.05 & 0 & 35.3\% & -2 & 745.1\% & -3 & 837.6\% & -1 & 224.3\% \\
GI AA & 60 & 18.25 & 0 & -1.6\% & 0 & 2.4\% & 0 & 8.5\% & 0 & -20.5\% \\
GI AB & 60 & 9.92 & 0 & -4.3\% & 0 & 6.4\% & 0 & 8.9\% & 0 & -20.2\% \\
GI BA & 59 & 71.88 & 0 & -13.9\% & 0 & -0.5\% & 0 & 1.0\% & 0 & -27.1\% \\
GI BB & 60 & 33.07 & 0 & -2.1\% & 0 & 4.2\% & 0 & 1.0\% & 0 & -19.2\% \\
Random & 3840 & 0.05 & 0 & 60.0\% & -1 & 300.0\% & -2 & 620.0\% & 0 & 0.0\% \\
Scholl & 1210 & 0.06 & 0 & 0.0\% & 0 & -16.7\% & 0 & 0.0\% & 0 & -16.7\% \\
Schwerin & 200 & 0.02 & 0 & 0.0\% & 0 & 0.0\% & 0 & 0.0\% & 0 & 0.0\% \\
Waescher & 17 & 0.31 & 0 & -6.5\% & 0 & -12.9\% & 0 & -9.7\% & 0 & -12.9\% \\
\bottomrule
\end{tabular}
\end{table}

\begin{table}[H]
\centering
\small
\caption{Comparison between the remaining versions and final version using time limit of 600s.}
\label{Table3b}
\centering
\begin{tabular}{lrrrrrrrrrrrrrr}
\toprule
& \multicolumn{2}{c}{\textsc{\mline{Final}{Version}}}  &\multicolumn{2}{c}{\textsc{\mline{Standard}{EPS}}}& \multicolumn{2}{c}{\textsc{\mline{No Dual Ineq.}{No B. Princing}}} & \multicolumn{2}{c}{\textsc{No MCRC}} &\multicolumn{2}{c}{\textsc{No Grouping}}\\
Name        &   Opt & Time    &$\Delta$ Opt & $\Delta$ Time&$\Delta$ Opt & $\Delta$ Time&$\Delta$ Opt & $\Delta$ Time&$\Delta$ Opt & $\Delta$ Time\\
\midrule
AI202 & 50 & 0.41 & 0 & -17.1\% & 0 & 4.9\% & 0 & 4.9\% & 0 & 36.6\% \\
AI403 & 50 & 6.01 & 0 & -24.3\% & 0 & -31.3\% & 0 & 5.7\% & 0 & 26.3\% \\
AI601 & 48 & 46.18 & 1 & 18.7\% & 1 & 16.5\% & 1 & -5.8\% & -3 & 61.8\% \\
AI802 & 47 & 93.67 & -2 & 33.7\% & -2 & 22.5\% & -3 & 20.2\% & -3 & 24.5\% \\
AI1003 & 38 & 190.60 & 1 & 0.5\% & -1 & 9.3\% & 1 & -4.2\% & 2 & -16.4\% \\
ANI201 & 50 & 0.48 & 0 & 0.0\% & 0 & -2.1\% & 0 & 12.5\% & 0 & 22.9\% \\
ANI402 & 50 & 2.69 & 0 & 4.1\% & 0 & 6.7\% & 0 & 18.6\% & 0 & -0.4\% \\
ANI600 & 50 & 13.42 & 0 & -5.6\% & 0 & 8.1\% & 0 & 42.0\% & 0 & -28.8\% \\
ANI801 & 49 & 38.78 & -1 & 35.3\% & 0 & 1.0\% & 0 & 13.7\% & 1 & -18.1\% \\
ANI1002 & 44 & 139.95 & -15 & 148.9\% & 0 & 4.0\% & -1 & 6.6\% & 1 & -18.7\% \\
FalkenauerT & 80 & 0.10 & 0 & 10.0\% & 0 & 10.0\% & 0 & 10.0\% & 0 & 1030.0\% \\
FalkenauerU & 80 & 0.02 & 0 & 0.0\% & 0 & 0.0\% & 0 & 0.0\% & 0 & 19350.0\% \\
Hard & 28 & 7.05 & 0 & 25.1\% & 0 & -6.1\% & 0 & 57.0\% & 0 & 46.7\% \\
GI AA & 60 & 18.25 & 0 & 14.3\% & 0 & 40.2\% & 0 & 10.7\% & -20 & 1726.2\% \\
GI AB & 60 & 9.92 & 0 & 9.0\% & 0 & 258.4\% & 0 & 10.4\% & -39 & 4748.7\% \\
GI BA & 59 & 71.88 & 0 & -12.3\% & 0 & 29.8\% & 0 & -3.1\% & -- & -- \\
GI BB & 60 & 33.07 & 0 & -7.2\% & 0 & 217.0\% & 0 & -6.4\% & -- & -- \\
Random & 3840 & 0.05 & 0 & 0.0\% & 0 & 20.0\% & 0 & 40.0\% & -2 & 3000.0\% \\
Scholl & 1210 & 0.06 & 0 & -16.7\% & 0 & 16.7\% & 0 & 0.0\% & -1 & 6383.3\% \\
Schwerin & 200 & 0.02 & 0 & 0.0\% & 0 & 50.0\% & 0 & 0.0\% & 0 & 1650.0\% \\
Waescher & 17 & 0.31 & 0 & -9.7\% & 0 & 54.8\% & 0 & -16.1\% & 0 & 1067.7\% \\
\bottomrule
\end{tabular}
\end{table}

\section{Detailed results for the CSP with 1 hour~\label{app::exp1h}}

In Table~\ref{Table4}, we present detailed results for the complete version using a time limit of one hour. We categorize the instances into those solved in the root node without RF heuristic (Opt\textsubscript{root}), those solved by executing RF heuristic in root node (Opt\textsubscript{RF}), those solved using up to 5 branches (Opt\textsubscript{$\leq$ 5}), and those solved using more than 5 branches (Opt\textsubscript{$>$ 5}). Additionally, we provide the average time spent in the pricer (Pricing Time), the average time Gurobi took to solve the RLM models (RLM Time), the average number of pricing calls (Pricing Calls), the average number of pricing calls that generated at least one column (Gen. Pricing Calls), and the average number of columns generated by each pricing call that generated at least one column (Cols by Gen. Pricing Calls).

Note that many instances are solved in the root node without RF heuristic, which indicates that the SCF formulation fortified by SR Cuts, BFD heuristic (used to produce the first incumbent solution), the rounding heuristic, and the smaller tolerance $\epsilon$ are enough to solve them. Moreover, our framework solves 180 of 250 AI instances and 198 of 250 ANI instances using up to 5 branches. We conjecture that the performance of our algorithm in these classes is due to the relaxation polyhedron for these instances being almost equal to the convex hull of the integer solution. This hypothesis is supported by the high Polyhedron Integrality~(PI) Ratio observed in these classes, especially in the ANI class, where many instances exhibit a PI Ratio around 80\%. This property justifies solving many instances using a few branches without the CRF heuristic.

\begin{table}[t!]
\centering
\caption{Detailed results for the final version using a time limit of one hour.}
\label{Table4}
\resizebox{\textwidth}{!}{%
\begin{tabular}{cccccccccccccc}
\toprule
Name & Total & Opt&  Opt\textsubscript{root} & Opt\textsubscript{RF} & Opt\textsubscript{$\leq$ 5} & Opt\textsubscript{$>$ 5}&  Time&  \mline{Pricing}{Time} & \mline{RLM}{Time} &  \mline{PI}{Ratio}  & \mline{Pricing}{Calls} &  \mline{Gen. Pricing}{Calls} & \mline{Cols by Gen.}{Pricing Call}\\
\midrule
AI202 & 50 & 50 & 12 & 26 & 3 & 9 & 0.4 & 0.05 & 0.3 & 48.3 & 50.0 & 38.1 &  26.3 \\
AI403 & 50 & 50 & 23 & 16 & 0 & 11 & 5.0 & 0.9 & 3.7 & 69.0 & 115.0 & 75.0 & 34.9 \\
AI601 & 50 & 50 & 18 & 12 & 1 & 19 & 57.1 & 9.2 & 43.0 & 52.9 & 378.8 & 180.7 & 26.4 \\
AI802 & 50 & 48 & 23 & 14 & 1 & 10 & 223.7 & 39.8 & 167.5 & 49.0 & 694.6 & 284.9 & 24.8 \\
AI1003 & 50 & 42 & 18 & 12 & 1 & 11 & 794.9 & 190.4 & 565.3 & 50.7 & 1234.8 & 487.3 & 20.5 \\
ANI201 & 50 & 50 & 31 & 0 & 9 & 10 & 0.4 & 0.07 & 0.3 & 77.7 & 60.0 & 46.8 & 22.1 \\
ANI402 & 50 & 50 & 39 & 0 & 6 & 5 & 2.2 & 0.6 & 1.5 & 86.3 & 69.0 & 57.1 & 42.7 \\
ANI600 & 50 & 50 & 28 & 0 & 9 & 13 & 11.8 & 3.6 & 7.6 & 85.4 & 115.6 & 86.1 & 47.9 \\
ANI801 & 50 & 50 & 35 & 0 & 7 & 8 & 57.0 & 15.0 & 39.5 & 84.4 & 192.0 & 122.8 & 50.1 \\
ANI1002 & 50 & 47 & 24 & 0 & 10 & 13 & 374.0 & 96.9 & 262.9 & 81.4 & 575.9 & 253.6 & 34.7 \\
FalkenauerT & 80 & 80 & 17 & 54 & 3 & 6 & 0.1 & 0.01 & 0.07 & 54.5 & 20.6 & 14.6 & 35.4 \\
FalkenauerU & 80 & 80 & 73 & 6 & 0 & 1 & 0.01 & 0.001 & 0.009 & 90.2 & 8.5 & 7.9 & 11.2 \\
Hard & 28 & 28 & 0 & 7 & 1 & 20 & 6.4 & 0.4 & 4.6 & 39.2 & 752.8 & 88.9 & 10.5 \\
GI AA & 60 & 60 & 38 & 22 & 0 & 0 & 14.5 & 14.0 & 0.4 & 94.5 & 19.6 & 15.8 & 80.7 \\
GI AB & 60 & 60 & 50 & 10 & 0 & 0 & 7.9 & 7.4 & 0.5 & 90.5 & 10.9 & 9.8 & 104.1 \\
GI BA & 60 & 59 & 37 & 22 & 0 & 0 & 102.0 & 100.5 & 1.3 & 94.4 & 45.9 &  16.5 & 76.5 \\
GI BB & 60 & 60 & 55 & 5 & 0 & 0 & 25.7 & 25.0 & 0.7 & 90.5 & 11.0 & 9.9 & 109.8 \\
Random & 3840 & 3840 & 3547 & 287 & 0 & 6 & 0.05 & 0.005 & 0.03 & 91.3 & 11.1 & 9.0 & 23.7 \\
Scholl & 1210 & 1210 & 1103 & 106 & 0 & 1 & 0.05 & 0.02 & 0.03 & 84.2 & 7.9 & 6.2 & 17.9 \\
Schwerin & 200 & 200 & 183 & 17 & 0 & 0 & 0.02 & 0.003 & 0.01 & 32.5 & 11.4 & 10.4 & 8.5 \\
Waescher & 17 & 17 & 4 & 9 & 1 & 3 & 0.3 & 0.1 & 0.1 & 9.4 & 67.2 & 58.7 & 6.6 \\
\bottomrule
\end{tabular}
}
\end{table}

Observe that our pricer is very efficient in AI and ANI classes since the most time spent is just for solving the RLM models. Moreover, GI instances have a very costly pricing problem, which is justified by the high number of items (around 5000 in the largest ones) and huge bin capacities, as ${W \in \{5\cdot 10^5, 1.5\cdot 10^6\}}$.

Finally, observe that our multiple-pattern generation usually adds dozens of patterns in each iteration. In fact, as we observed, the number of patterns generated by the complete version and the best pattern version is similar, thus these multiple patterns generated are really reducing the number of pricing iterations.

\section{Our benchmark\label{app::benchmark}}

Excluding around 30\% of AI instances, the BPP Lib mainly consists of instances that our algorithm can readily solve, even when employing only our Rounding Heuristic. Therefore, we are confident that employing strong relaxation formulations alongside any reasonable branching scheme embedded in a node processing rule that induces a diving heuristic is sufficient to find an optimal solution for these instances. However, given the strong NP-hardness of the problem, finding an optimal solution for large instances should be a challenge greater than to executing simple variable fixing heuristics a few times.

For variable fixing heuristics, a hard IRUP instance is characterized by the ease with which primal heuristics select a subset of patterns for fixing in a way that residual instances are non-IRUP, and by the difficulty of identifying the residual instances as non-IRUP\@. Although this subset of patterns cannot be part of an optimal solution, the variable fixing heuristics must select them, thus, these variables should have positive values in the RLM solution. Moreover, we posit that these non-IRUP residual instances are more prevalent in cases with a highly fractional polyhedron, where optimal RLM solutions typically comprise few integer positive variables and numerous positive variables. In the following, we outline our approach to generating instances that meet these criteria.

Given integer numbers $n$ and $W$, where $n = 3l \cdot 3 ^k$, and $l, k \in \mathbb{Z}_+^*$ with $3 \leq l < 9$, we construct a set $S$ consisting of $l$ triples of items with a total size sum equal to $W$. The size $w_i$ of the first item in each triple is uniformly at random chosen from the range $[\frac{W}{5}, \frac{2W}{5}]$. The size of the second item is chosen from $[\frac{W}{5}, W - w_i - \frac{2W}{5}]$, and the size of the last item is set to the remaining capacity. Additionally, to create more challenging instances, if any new item in a triple is identical to any existing item in the instance, it is discarded, and another attempt is made. The current triple is accepted only on the $100^{\text{th}}$ try to avoid looping.

Subsequently, we iterate the following process $k$ times, using the previous set $S$ as a starting point to generate a new set $S'$. For each triple of items $(i_1, i_2, i_3) \in S$, we divide these items into three new triples following the same procedure, except that the first item of each triple is set to $i_j$. After all iterations, we obtain $n$ items, and an optimal solution consists solely of full bins.

Note that for an item $i$ created in iteration $r$, there are at least ${k - r + 1}$ complete patterns in $\mathcal{P}$ that contains $i$, contributing to the creation of a highly fractional polyhedron. Additionally, our initial experiments reveal that both the bin capacity $W$ and the range of item sizes also influence the integrality of the relaxation polyhedron. With this in mind, we generated $50$ instances for each ${(n, W) \in {(216, 10^3), (405, 1.5 \cdot 10^3), (648, 2 \cdot 10^3)}}$.

\begin{table}[b!]
\centering
\small
\caption{Results for our CSP benchmark using our framework (time limit of 3600s).}
\label{Table5}
\resizebox{\textwidth}{!}{%
\begin{tabular}{cccccccccccccc}
\toprule
N  & Total & Opt &  $\text{Opt}_{>5}$&  Time &  \mline{Pricing}{Time} & \mline{RLM}{Time}  & Cols    &    Cuts  &PI ratio &     \mline{Pricing}{Calls} &     \mline{Gen. Pricing}{Calls} & \mline{Cols by Gen.}{Pricing Call}\\
\midrule
216 & 50 & 44 & 44 & 1082.4 & 220.2 & 730.6 & 3251.3 & 554.6 & 4.6 & 48891.0 & 1119.0 & 2.9 \\
405 & 50 & 34 & 34 & 2026.2 & 552.0 & 1225.9 & 12126.7 & 897.1 & 2.7 & 37241.9 & 3458.6 & 3.5 \\
648 & 50 & 13 & 13 & 3207.3 & 925.5 & 1893.4 & 20078.6 & 1273.6 & 2.9 & 28314.8 & 4372.2 & 4.6 \\
\bottomrule
\end{tabular}
}
\end{table}

Table~\ref{Table5} presents the results obtained by our framework\footnote{
    Unfortunately, we do not have access to the NF-F algorithm of~\citet{Loti_2022} to compare our results.
    }
in this new benchmark. We observed that the PI ratio is, on average, less than 5\% in the root node of all generated instances, directly impacting the number of solved instances, as we anticipated. All solved instances required more than five branches to find the optimal solution, with the most costly step being solving RLM models. Furthermore, given the high discrepancy between the number of pricing calls and the number of pricing calls that generated at least one column, indicating that the B\&B nodes generally require few or no new patterns to be solved. This is also supported by the number of columns generated by pricing calls that generated at least one column, which is relatively small compared to other benchmarks.

Finally, even using a time limit of one hour, our framework found an optimal solution for 36 of 50 instances with 216 items and left many instances open with 405 and 648 items. Therefore, these instances are a new and challenging CSP benchmark that can be considered in future works and can probably be adapted to other problems.

\section{Binary Search Algorithm from IPMS Solver\label{app::IPMS}}
Next, we present in Algorithm~\ref{alg::binarysearch} the pseudocode from our IPMS Solver, which consists of our CSP solver embedded in a binary search.  For the sake of simplicity, we use the notation for CSP, representing the set of jobs $J_{j}$ with processing time $s_j$ as an item $i \in I$ with size $w_i = s_i$ and demand $d_i =   \vert J_{j} \vert $.

\begin{algorithm}[ht]
    \small
    \SetKwInOut{Input}{input}\SetKwInOut{Output}{output}
    $LB = \max(\lceil \sum_{i \in I} \frac{d_i w_i}{m} \rceil, \max_{i \in I} w_i)$ \hspace{0.5cm} \tcp{Lower bound}

    $inc \gets$ LongestProcessingTime$(I, m)$ \hspace{0.5cm} \tcp{ Incumbent solution}

    $UB \gets getValue(inc)$ \hspace{0.5cm}\tcp{Upper bound}

    $p \gets 0$

    $isFirstStage \gets True$

    \While{$LB \neq UB$}{
        \eIf{$isFirstStage$}{
            $W \gets  \min(LB + 2^P - 1, UB - 1)$

            $p \gets p + 1$
        }{
            $W = \lfloor \frac{LB + UB}{2} \rfloor$
        }

        $\Zi \gets$ CSPSolver~($I, W$)

        \eIf{$z(\Zi) \leq m$}{

            $UB \gets W$

            $inc \gets \Zi$

            $isFirstStage \gets False$

        }{
            $LB \gets W + 1$
        }
    }

    \textbf{return} inc
  \caption{\label{alg::binarysearch}IPMS Solver}
\end{algorithm}

\end{appendices}

\end{document}